\documentclass[12pt]{article}

\usepackage{amsmath,amssymb,amsthm}

\newtheorem{theorem}{Theorem}[section]
\newtheorem{lemma}{Lemma}[section]
\newtheorem{proposition}{Proposition}[section]
\newtheorem{corollary}{Corollary}[section]

\theoremstyle{definition}
\newtheorem{example}{Example}[section]
\newtheorem{remark}{Remark}[section]
\newtheorem{definition}{Definition}[section]

\newcommand{\bA}{\mathbf{A}}
\newcommand{\bF}{\mathbf{F}}
\newcommand{\bH}{\mathbf{H}}
\newcommand{\bN}{\mathbf{N}}
\newcommand{\bS}{\mathbf{S}}
\newcommand{\bU}{\mathbf{U}}
\newcommand{\bZ}{\mathbf{Z}}

\newcommand{\bC}{\mathbb{C}}
\newcommand{\bR}{\mathbb{R}}

\newcommand{\cA}{\mathcal{A}}
\newcommand{\cC}{\mathcal{C}}
\newcommand{\cF}{\mathcal{F}}
\newcommand{\cH}{\mathcal{H}}

\newcommand{\cL}{\mathcal{L}}
\newcommand{\cM}{\mathcal{M}}
\newcommand{\cR}{\mathcal{R}}

\newcommand{\nn}{\nonumber}
\newcommand{\cplus}{\overset{\circ}{+}}
\newcommand{\hb}{\hbar}
\newcommand{\pa}{\partial}
\newcommand{\wh}{\widehat}
\newcommand{\ve}{\varepsilon}
\newcommand{\ol}{\overline}
\newcommand{\oz}{\overline{z}}
\newcommand{\uF}{\underline{F}}
\newcommand{\opa}{\overline{\partial}}
\newcommand{\od}{\overset{\rm def}{=}}
\newcommand{\cPhi}{\overset{\circ}{\Phi}}

\renewcommand{\Re}{\operatorname{Re}}
\renewcommand{\Im}{\operatorname{Im}}

\newcommand{\ad}{\operatorname{ad}}
\newcommand{\su}{\operatorname{su}}
\newcommand{\so}{\operatorname{so}}
\newcommand{\const}{\operatorname{const}}

\newcommand{\fS}{\operatornamewithlimits{\mathfrak S}}

\begin{document}

\title{Noncommutative Algebras, Nano-Structures, 
and Quantum Dynamics Generated by Resonances}

\author{Mikhail Karasev\thanks{This work was partially supported
by RFBR (grant 05-01-00918-a) and by INTAS (grant 00-257).}\\ \\
\it Moscow Institute of Electronics and Mathematics\\
karasev@miem.edu.ru}

\date{}
                                                
\maketitle

\begin{abstract}
We observe ``quantum'' properties of resonance
equilibrium points and resonance univariant submanifolds in the
phase space. Resonances between Birkhoff or Floquet--Lyapunov
frequencies generate quantum algebras with polynomial
commutation relations. Irreducible representations and
coherent states of these algebras correspond to certain
quantum nano-structure near the classical resonance motion.
Based on this representation theory and nano-geometry,  
for equations of Schr\"odinger or wave type in various
regimes and zones (up to quantum chaos borders) we describe the
resonance spectral and long-time asymptotics, resonance
localization and focusing, resonance adiabatic and spin-like
effects. 
We discuss how the mathematical phase space nano-structures
relate to physical nanoscale objects like dots, quantum
wires, etc.
We also demonstrate that even in physically macroscale Helmholtz
channels the resonance implies a specific quantum character of
classical wave propagation. 
\end{abstract}

\tableofcontents

\bigskip

The paper consists of several parts. 
Part~I follows the material of 
the author lectures 
at Petrovskii seminar \& Moscow Math. Society Conference
(May, 2004) 
and St.-Peters\-burg University \& Steklov Math. Institute
Conference  
(June, 2004).

Part~II contains the systematic description of resonance
algebras. 
In the next parts, 
we shall analyze the nano- and micro- phase space structures 
and discuss various types resonance phenomena for the
Schr\"odinger type  and  wave equations, 
and also establish a bridge to real physical scales. 

\part{}

\setcounter{section}{-1}

\section{Introduction}%0

Equations of mathematical physics describing propagation of
waves admit solutions (or regimes, or states) of very different
kinds. Many important physical applications deal with 
solutions which are not completely chaotic 
but follow distinguished ``integrable'' motions, say, 
equilibria.
These solutions can be considered as certain
excitations around the integrable  classical core motion. 
The leading part of excitations is described by a model
equation. 
By resolving this model equation, it is then possible to compute
solutions of the original problem via the perturbation theory. 
Such a scheme goes back to Laplace, Rayleigh, Poincare, Ehrenfest,
Birkhoff, Bogolyubov.

Usually, in this approach one presupposes to obtain 
the model equation to be as simple as possible, 
say, reducible to trivial scalar operators 
or to first-order differential operators of classical type.  
But, in many crucial cases, the model equation occurs to be  
of a nontrivial quantum type, that is, it 
carries a certain noncommutative algebra structure. 
The most typical reason for these quantum algebras to appear 
is the degeneracy of spectrum of the core integrable motion.

Such a quantum behavior can arise not only in
nanoscale problems of atomic physics, but even in purely
classical wave propagation problems at usual macrophysical
scales. Thus one can claim that there are quantum effects 
in nonquantum wave systems as well. This very interesting
phenomenon is the main motive of our present work.

Studying wave equations very often uses 
the analogy with classical mechanical systems and exploits 
the phase space geometry. The quantization technique, 
developed to obtain an operator representation of the phase
space geometric structures, 
can be effectively applied to construct approximate
or exact solutions of wave equations. 
This general claim is supported, in particular, 
by the progress in the semiclassical asymptotics~\cite{GuSt}.

The ray method \cite{kar1} and the general Maslov's canonical
operator theory~\cite{kar2} were developed to construct
asymptotic (semiclassical) solutions of PDE localized at points,
trajectories, tori, or other invariant submanifolds in a phase
space.  

In multidimensional case all approaches 
of the semiclassical approximation theory 
have 
the well-known stumbling block: 
the incommensurability condition for frequencies of the
classical Hamiltonian dynamics over the invariant submanifolds.

If the resonance between frequencies takes place 
then the usual methods fail. The resonance situation is an old
open problem in the theory of wave and quantum equations. 
Of course, this is related to the
resonance problems in classical mechanics~\cite{kar14},
for instance, to the resonance theory of averaging and
normal forms~\cite{kar14,kar15}.
However, in the wave mechanics, 
at least in stable case, the resonance problems are simpler,
since they deal with a discrete spectrum and a finite-dimensional
degeneracy.

An example can be presented by the 
Schr\"odinger (or Helmholtz) operator 
whose potential (or index of refraction)
has a nondegenerate minimum (maximum) point 
and square roots of eigenvalues of the second derivative matrix 
at this point, in Euclidean coordinates, are commensurable. 
The semiclassical asymptotics of the spectrum near
such a bottom (or maximum) is an intriguing question
unsolved until now.

Another known unsolved problem: 
spectral asymptotics corresponding 
to stable trajectories, e.g., geodesics, 
in the case of a resonance between
Lyapunov frequencies or 
Floquet frequencies~\cite{kar1,kar16}. 

One more interesting question is the long-time evolution of wave
packets localized at a resonance stable equilibrium point, 
or trajectory, or torus. 

We suggest a way to solve these resonance problems by
studying some noncommutative algebraic structures 
in micro- and nano-zones near the resonance core motion, 
and by applying and developing general methods of quantum
geometry.
We observe that each resonance proportion 
between frequencies generates an algebra with 
polynomial commutation relations 
(polynomial Poisson tensor) 
and its irreducible representations 
are given by hypergeometric K\"ahlerian 
structures. The model equation is an equation 
over this resonance algebra. 
It is not of classical type, i.e., 
it has an order greater than one 
(in the irreducible representation). 
Therefore the nanozone becomes a purely quantum one. 
In microzones the model equation is reduced 
to the first order 
and resolved by the semiclassical technique.
We discuss some important physical examples 
and detect interesting effects generated 
by resonances in classical and quantum wave equations.

The given first part of the paper describes mostly the material
of the author's  lectures \cite{kar3,kar21}.

\section{Correlation of modes in resonance clusters}%1

Let us consider the operator 
\begin{equation}\label{k1}                       
\bH=\bH_0+v\qquad\text{in}\quad L^2(\bR^2),
\end{equation}
where 
\begin{equation}\label{k2}                      
\bH_0=-\frac{\hb^2}{2}\Delta+V_2,\qquad v=V_3+V_4+\dots,
\end{equation}
and by $V_j$ we denote a $j$-linear form on $\bR^2$.
For instance, $V_2$ is a quadratic form, and the operator
$\bH_0$ is just the Hamiltonian of a harmonic oscillator.

Let $A>0$ and assume that the domain in $\bR^2$ 
where the potential $V_2+v$ does not exceed the value $A$ is a
connected neighborhood of zero. 
We are interested in spectral properties of the operator
$\bH$ on the energy interval $(0,A)$.

Denote by $\alpha,\beta$ the frequencies of the
oscillator~$\bH_0$.  
Then the spectrum of $\bH_0$ consists of the numbers
\begin{equation}\label{k3}                      
\lambda_{m,l}=\hb\alpha\Big(m+\frac12\Big)
+\hb\beta\Big(l+\frac12\Big)
\end{equation}
with known eigenfunctions $| m,l\rangle\in L^2(\bR^2)$ given by
Hermite polynomials multiplied by the Gaussian exponent.

The resonance occurs when the following condition holds:
\begin{equation}\label{k4}                      
\frac{\alpha}{\beta}=\text{rational number}.
\end{equation}
In this case, the spectrum of $\bH_0$ is degenerate, that is,
there are nontrivial clusters of pairs $(m,l)\sim(m',l')$ 
such that $\lambda_{m,l}=\lambda_{m',l'}$.

To understand what is the effect of resonance, let us consider
the matrix elements
$$
E(t)=\langle m',l' | e^{-\frac{it}{\hb}\bH}| m,l\rangle
$$
which control the correlation (the transition probability) 
between the modes $| m,l\rangle$ and $| m',l'\rangle$. 

In the nonresonance case where the ratio $\alpha/\beta$ is
irrational, for any $t\sim\hb^{-s}$, one has
$$
E(t)\sim \hb^{1/2}\qquad\text{as}\quad \hb\to0.
$$
But in the resonance case, if $t\sim \hb^{-s}$ (for some~$s$) 
and the pairs $(m,l)\sim(m',l')$
are inside a cluster, then 
$$
E(t)\sim O(1)\qquad\text{as}\quad \hb\to0.
$$

{\it Thus the resonance implies a strong correlation of modes inside
the cluster, or a possibility of transition along modes in the
cluster, that is, a dynamics.
A phase space geometry which underlies this dynamics is of
compact type, since the cluster is finite.}

These are analytic and geometric consequences of the resonance.
But behind all there is an algebraic phenomenon generated by the
resonance.

\section{Noncommutative resonance algebras}%2

Let us consider a neighborhood of the origin in $\bR^2$ 
of order $\hb^{1/N}$. 
We call this domain an {\it $N$th microzone\/} if $N>2$.
If $N=2$, then we use a specific term a {\it nanozone}.

In micro or nanozones, the potential $v$ in (\ref{k1}), (\ref{k2}) 
can be considered as a perturbation with respect to the leading
part $\bH_0$. 

Denote by $M_0$ the algebra of integrals of motion for
$\bH_0$, or the commutant of $\bH_0$. Thus, each element from
$M_0$ commutes with $\bH_0$. 
In the $N$th microzone, 
one can find a unitary operator $\bU$ such that 
\begin{equation}\label{k5}                      
\bH=\hb^{2/N}\bU^{-1}
(\bH_0+\hb^{1/N}\bF_1+\hb^{2/N}\bF_2+\dots+\hb^L\bF_{LN})\bU
+O(\hb^{L+3/N}),
\end{equation}
where $\bF_j\in M_0$ for all $j$. 
For instance, the operator $\bF_1$ is just the projection 
of $V_3$ onto $M_0$, the operator $\bF_2$ is determined 
(explicitly) by $V_3$ and $V_4$, etc.
The representation (\ref{k5}) is the result of application of the
quantum averaging method developed
in a similar framework in \cite{We,kar4,kar5,kar22}. 
As we see from (\ref{k5}),
the study of the operator $\bH$ in the $N$th microzone can be 
reduced to the study of operators from the algebra~$M_0$.

In the nonresonance case the algebra $M_0$ is commutative.

{\it In the resonance case} (\ref{k4}), 
{\it the algebra $M_0$ is noncommutative}. 

Thus, the resonance implies the ``quantum'' behavior of the
problem under study. The word ``quantum'' we use as a synonym 
of ``noncommutative.''
The ``quantum ray method,'' the ``quantum characteristics,''
the ``quantum geometry''~-- all this appears as a consequence 
of noncommutativity of the algebra~$M_0$ under the resonance.

The commutation relations in the algebra $M_0$ are of the
following type:
\begin{equation}\label{k6}                      
[\bA_j,\bA_k]=-i\hb'\Psi_{jk}(\bA).
\end{equation}
Here $\hb'=\hb^{1-2/N}$, and $\bA=(\!(\bA_j)\!)$ 
is a finite set of generators, 
$\Psi_{jk}$ is a Poisson (quantum) tensor. One can
choose the generators in such a way that the components
$\Psi_{js}$ be polynomial, and so we can say that $M_0$ 
is an algebra with polynomial commutation relations. In general,
these relations do not belong to the class of Lie algebras and
present more complicated finitely generated algebras 
whose study began not so long ago (see the reviews in
\cite{kar6,kar7} and the references therein).

In the case of the simplest isotropic resonance $1:1$
(that is, $\alpha=\beta=1$), one has the following relations 
in~$M_0$: 
\begin{align}
[\bA_1,\bA_2]&=-i\hb'\bA_3,
\nonumber\\
[\bA_2,\bA_3]&=-i\hb'\bA_1,
\tag{6a}\\
[\bA_3,\bA_1]&=-i\hb'\bA_2.
\nonumber
\end{align}
So, in this case, $M_0$ is the enveloping of the Lie algebra
$\su(2)$.  
The spectral analysis of the operator (\ref{k5}) 
in this case is reduced to the study of a Hamiltonian over 
$\su(2)$,
which can be done by the standard technique 
(see, for instance, \cite{per,kurc,kar23}).

Much more interesting case is represented by anisotropic 
resonances, say, the resonance $1:2$
(where $\alpha=1$, $\beta=2$). Denote by $x',y'$ 
rescaled Cartesian coordinates adapted to the $N$th microzone. 
Introduce the annihilation operators 
$$
\eta=x'+\hb'\frac{\partial}{\partial x'},\qquad
\zeta=y'+\hb'\frac{\partial}{\partial y'}.
$$
Then the algebra $M_0$ is generated by self-adjoint operators 
\begin{align}
\bA_1&=\frac14\eta^*\eta,\qquad 
\bA_2=\frac1{12}(\eta^*\eta-4\zeta^*\zeta),
\label{k7}\\
\bA_3&=\frac18(\zeta^*\eta^2+{\eta^*}^2\zeta),\qquad 
\bA_2=\frac1{8i}(\zeta^*\eta^2-{\eta^*}^2\zeta),
\nonumber
\end{align}
The commutation relations (\ref{k6}) in this case 
(resonance $1:2$) are the following:
\begin{align}
[\bA_1,\bA_2]&=0,\qquad [\bA_1,\bA_3]=-i\hb'\bA_4,\qquad
[\bA_1,\bA_4]=i\hb'\bA_3,\nonumber
\\
[\bA_2,\bA_3]&=-i\hb'\bA_4,\qquad [\bA_2,\bA_4]=i\hb'\bA_3,
\label{k8}                      
\\
[\bA_3,\bA_4]&=-3i\hb'\bigg(\bA_1\bA_2-\frac{\hb'}{4}\bA_1
+\frac{\hb'}{4}\bA_2\bigg)
\nonumber
\end{align}
The Casimir elements of this non-Lie algebra are
\begin{equation}
\begin{aligned}
\bC_1&=\bA_1-\bA_2,
\\
\bC_2&=3\bA_1^2 \bA_2  -  \bA_1^3  +  \bA_3^2  
+ \bA_4^2 - \frac{3 \hbar'}{2} \bA_1^2  
+ \frac{3 \hbar'}{2} \bA_1 \bA_2  
+ \frac{3 \hbar^{\prime 2}}{4} \bA_2   
+ \frac{\hbar^{\prime 2}}{4} \bA_1.   
\end{aligned}
\tag{8a}
\end{equation}

Relations (\ref{k8}) remind quadratic algebras appearing 
in the theory of infinite-dimensional integrable systems
(but now without the Hopf axiom, see~\cite{Fadd,kar8}).

{\it Each resonance proportion $q:r$ in {\rm(\ref{k4})} 
generates an algebra $M_0$ with polynomial commutation relations
{\rm(\ref{k6})}.
The coefficients of the polynomial $\Psi_{jk}$ 
in {\rm(\ref{k6})} can be made all be integer numbers 
determined by $q,r$. 
The same is true for more than two resonance frequencies as well.}

From (\ref{k5}) it follows that, in the micro- and nanozone near the
bottom of the potential, the operator (\ref{k1}) is approximately
reduced to the Hamiltonian 
\begin{equation}\label{k9}
\hb^{2/N}\big(\bH_0+\hb^{1/N}f_1(\bA)+\hb^{2/N}f_2(\bA)+\dots\big),
\end{equation}
with some polynomials $f_1,f_2,\dots$ in the generators
$\bA$ of the algebra $M_0$. 
The leading part $\bH_0$ in (\ref{k9}) is a Casimir element in
$M_0$ and can be replaced by a scalar (in the irreducible
representation). 

For example, consider the case of the resonance $1:2$. 
The operator
\begin{equation}\label{k10}                      
\bH=-\frac{\hb^2}{2}\Delta+\frac12x^2+2y^2+x^2y
+\gamma x^4,\qquad \gamma\geq \frac18,
\end{equation}
in the $N$th nanozone is reduced to 
\begin{equation}\label{k11}                      
\hb\Big(n+\frac32\Big)+\frac{\hb^{3/N}}{\sqrt{2}}\bA_3+O(\hb^{4/N}),
\end{equation}
where $\bA_3$ is the generator of the algebra (\ref{k8})
given by the third formula in~(\ref{k7}).
Thus, one needs to study the spectrum of the element $\bA_3$ in
the $n$th irreducible representation of the algebra (\ref{k8}) 
in order to find the bottom part of the spectrum of the
operator~(\ref{k10}). 

The following important question arises: how to construct
irreducible representations of algebras of the type~(\ref{k6})?

For the case of Lie algebras
(in particular, for the resonance case $1:1$), 
one has the Kirillov orbit method \cite{Kir} and 
the general Kostant--Souriau geometric quantization
\cite{Kos,Sou}. 
For anisotropic resonances and algebras
with nonlinear commutation relations like~(\ref{k6}), 
one needs an extension of this quantization scheme.

\section{Quantum nano-geometry}%3

The directing idea of the geometric quantization is to use the
classical Poisson (symplectic) geometry and polarizations 
in order to determine quantum objects. The classical analog 
of the algebra $M_0$ is a Poisson algebra $\cM_0$ of functions
on the phase space commuting with the symbol of the
operator~$\bH_0$. The set $\cM_0$, of course, was considered in
classical mechanics dealing with resonance systems, see,
e.g.,~\cite{Gust}, but the Poisson structure on $\cM_0$ 
was not described and studied.

The classical analog of the quantum commutation relations
(\ref{k6}) is given by a Poisson tensor of a
polynomial type.
For example, in the case of resonance $1:2$, 
we obtain the following quadratic Poisson algebra:
\begin{align}
\{\cA_1,\cA_2\}&=0,\qquad \{\cA_1,\cA_3\}=\cA_4,\qquad
\{\cA_1,\cA_4\}=-\cA_3,\nonumber
\\
\{\cA_2,\cA_3\}&=\cA_4,\qquad \{\cA_2,\cA_4\}=-\cA_3,
\label{k12}                      
\\
\{\cA_3,\cA_4\}&=3\cA_1\cA_2.
\nonumber
\end{align}
The Casimir functions in this algebra are  
$$
\cC_1=\cA_1-\cA_2, \qquad
\cC_2=3\cA^2_1\cA_2-\cA^3_1+\cA^2_3+\cA^2_4.
$$
Generic symplectic leaves are two-dimensional surfaces
\begin{equation}\label{k13}                      
\Omega=\{\cC_1=\text{const},\cC_2=\text{const}\}\subset \bR^4
\end{equation}
diffeomorphic to $\bS^2$. 
The leaves are K\"ahlerian manifolds
with respect to the complex structure 
\begin{equation}\label{k14}
z=\frac{\cA_3+i\cA_4}{c-\cA_1},\qquad c=\text{const}.
\end{equation}
The symplectic form $\omega_0$ generated by brackets
(\ref{k12}) on the leaf (\ref{k13}) can be expressed by a
K\"ahlerian potential $F_0$ in the standard way
\begin{equation}\label{k15}                      
\omega_0=i\ol{\partial}\partial F_0=ig_0d\ol{z}\wedge dz.
\end{equation}
Note that the differential $2$-form 
\begin{equation}\label{k16}                      
\rho=i\ol{\partial}\partial\log g_0
\end{equation}
is the Ricci form on $\Omega$ corresponding to the complex
structure (\ref{k14}). 

Now if one follows the geometric quantization ideas, the
linear bundle over $\Omega$ with the curvature $i\omega_0$ must
be introduced. 
In the Hilbert space of sections of this bundle, the operators
of irreducible representation would act.

However, there are two principal problems.
First, we do not know the measure on $\Omega$ with respect to
which the Hilbert norm in the space of sections has to be
defined. This measure must satisfy a 
{\it reproducing property}~\cite{kar9}.  
For the inhomogeneous case (where relations (\ref{k6}) are not
linear) the existence of such a reproducing measure is, 
in general, unknown.

Secondly, even if one knows the reproducing measure, the problem
is that the operators of irreducible representation of the
algebra (\ref{k6}) constructed canonically by the geometric
quantization scheme would be pseudodifferential, 
but not differential operators.

That is why we modify the quantization scheme and from the very
beginning replace the symplectic form $\omega_0$ by another
``quantum'' form $\omega$ in a way that guarantees the
existence of the reproducing measure and the existence
of irreducible representations of the algebra (\ref{k6})
by differential operators.
This approach is explained in~\cite{kar9}.

Note that the opportunity to obtain irreducible representations 
of the algebra $M_0$ by differential operators is exactly 
the reason why the polynomial structure of the right-hand side  
of~(6) is so critical.

\section{Irreducible representations and\\ 
coherent states of resonance algebra}%4

Here we demonstrate calculations for algebra (\ref{k8}) related
to the resonance $1:2$.

Let us consider the following hypergeometric equation 
\begin{align*}
&2\hb's\frac{d^2 K}{ds^2}+(s+\ve_n\hb')\frac{dK}{ds}
-[{n}/{2}]K=0,\qquad s>0
\\
&K(0)=1,
\end{align*} 
where $n\geq0$ is an integer, $\ve_n=1$ or $\ve_n=3$ 
if $n$ is even or odd, and 
the brackets $[\cdots]$ denote the integer part.
The solution is given by a hypergeometric function 
of ${}_1\cF_1$ type, 
or more precisely, 
\begin{equation}\label{k17}                      
K(s)=\sum^{[{n}/2]}_{j=0}
\frac{[{n}/2]!}{j! (2j-2+\ve_n)!! ([{n}/2]-j)!}
\bigg(\frac{s}{\hbar'}\bigg)^j.
\end{equation} 
Here the double factorial $!!$ denotes the product over odd 
numbers, starting from~$1$.    

Also consider the ``dual'' hypergeometric equation
$$
2\hb's\frac{d^2L}{ds^2}
-(s +\ve_n \hbar'-4\hbar') \frac{dL}{ds}
-\big([{n}/2]+2\big) L =0,\qquad s>0,
\leqno{\rm(17a)}
$$
$$
\frac1{\hbar'}
\int^\infty_0 L(s)\,ds = 1.
$$
Introduce the Hilbert space $\cH_n$ of polynomials 
square integrable over $\bC$ with respect to the measure
$$
d\mu=L(|z|^2)\, d\ol{z}\,dz,\qquad z\in\bC.
$$

\begin{lemma}%1
The operators 
\begin{align}
\check{A}_1&=\frac{\ve_n-1}{4} \hb'+\hb'\ol{z}\ol{\partial},
\qquad\text{where}\quad\ol{\partial}=\partial/\partial\ol{z},
\nonumber\\
\check{A}_2&=\frac{\ve_n-1}{12} \hb'
-\frac{2\hb'[{n}/{2}]}{3}
+\frac{5\hb'}{12}\ol{z}\ol{\partial},
\label{k18}\\
\check{A}_3&={\hb'}^2\ol{z}\ol{\partial}^2
-\frac{\hb'}{2}(\ol{z}^2-\ve_n\hb')\ol{\partial}
+\frac{\hb'[{n}/{2}]}{2}\ol{z},
\nonumber\\
\check{A}_4&=-i{\hb'}^2\ol{z}\ol{\partial}^2
-\frac{i\hb'}{2}(\ol{z}^2+\ve_n\hb')\ol{\partial}
+\frac{i\hb'[{n}/{2}]}{2}\ol{z}
\nonumber
\end{align}    
realize the Hermitian irreducible representation of the algebra
{\rm(\ref{k8})} in the Hilbert space~$\cH_n$.
In this representation the Casimir elements {\rm(8a)} are
$\check{C}_1={n\hbar'}/{3}$ and 
$\check{C}_2=0$.
\end{lemma}

\begin{lemma}%2
Let us define the ``vacuum'' vector $|0\rangle$ 
as the solution of the equations 
\begin{align*}
\bA_1 |0\rangle=a_1|0\rangle,\qquad
\bA_2 |0\rangle=a_2|0\rangle,\qquad
(\bA_3+i\bA_4) |0\rangle=0,
\end{align*}
where $a_1, a_2$ are constants given by $a_j=\check{A}_j 1$. 
Also define the coherent states in $L^2(\bR^2)$: 
\begin{equation}\label{k19}                      
|z\rangle\od K\big(z(\bA_3-i\bA_4)\big)|0\rangle.
\end{equation}    
Then the integral transformation 
\begin{equation}\label{k20}
T(\varphi)=\frac{1}{2\pi\hbar'}
\int_{\bC}\varphi(\ol{z})\,|z\rangle\,d\mu(\ol{z},z),
\qquad 
T:\, \cH_n\to L^2(\bR^2),
\end{equation}
intertwines the representation {\rm(\ref{k7})} of the algebra
{\rm(\ref{k8})} and the irreducible representation
{\rm(\ref{k18})}, i.e.,  
$$
\bA_j\circ T=T\circ \check{A}_j\qquad (j=1,\dots,4).
$$
\end{lemma}

\begin{lemma}%3
The hypergeometric polynomial $K$ {\rm(\ref{k17})} is 
the reproducing kernel for the space~$\cH_n$.
Moreover,
\begin{align*}
K(|z|^2)=\langle z|z\rangle,\qquad 
\frac{1}{2\pi\hbar'}\int_{\bC} |z\rangle \langle z|\,d\mu=P_n,
\end{align*}
where $P_n$ is the projection in $L^2(\bR^2)$ onto the $n$th
irreducible component.
\end{lemma}

\begin{lemma}%4
Let the quantum K\"ahlerian form on $\Omega$ be defined by 
$$
\omega=i\hb'\ol{\partial}\partial \ln K,
$$
where $K$ is the hypergeometric polynomial {\rm(\ref{k17})},
and $\hbar'=\hbar^{1-2/N}$, $N\geq 2$.
Then the reproducing measure on $\Omega$ is given by 
$$
dm=K(|z|^2)L(|z|^2)\,d\ol{z}\,dz,
$$
where $L$ is the solution of the dual hypergeometric
equation~{\rm(17a)}.

In a micro-zone, where $N>2$ and $\hbar'$ is a small parameter,  
$n\sim 1/\hbar'$,
the following asymptotics hold{\rm:}
$$
\omega=\omega_0+\frac{\hb'}{2}\rho+O({\hb'}^2),\qquad 
dm=dm_0(1+O(\hb')).
$$
Here $dm_0=g_0 d\ol{z}\,dz$ is the Liouville measure
corresponding to the form $\omega_0$ on $\Omega$,
and $\rho$ is the Ricci form, see {\rm(\ref{k15}), (\ref{k16})}.
Moreover, one has
$$
K=e^{F_0/\hb'}\sqrt{g_0}(1+O(\hb')),\qquad
L=e^{-F_0/\hb'}\sqrt{g_0}(1+O(\hb')).
$$

Also the following identities hold{\rm:}
\begin{equation}\label{k21}
\frac1{2\pi \hb'}\int_{\Omega}\omega=[{n}/2],
\qquad
\frac1{2\pi \hb'}\int_{\Omega} dm=[{n}/2]+1.
\end{equation}
\end{lemma}

From this lemma one can clearly see that 
there is an essential difference
between the classical and quantum K\"ahlerian structures on
symplectic leaves~$\Omega$.
In micro-zones, where $N>2$, this difference is just
asymptotical: 
the quantum structure is an $\hbar'$-perturbation of the
classical one. But in the nanozone, where $N=2$ and $\hbar'=1$,
the quantum structure is not a perturbation of the classical
structure. Thus, one can talk about a specific 
{\it quantum nano-geometry\/} accompanying the resonance.

\section{Spectrum asymptotics in resonance clusters}%5

Now coming back to the spectral problem for the Hamiltonian
(\ref{k10}), we can apply the coherent transform~(\ref{k20}) 
to the operator (\ref{k11}). Then our problem in the $N$th micro zone 
is reduced to studying the operator 
\begin{equation}\label{k22}
\hb\bigg(n+\frac32\bigg)+\frac{\hb^{3/N}}{\sqrt{2}}\check{A}_3
+O(\hb^{4/N}),
\end{equation}
where $\check{A}_3$ is given by (\ref{k18}).

Recall that $\hb'=\hb^{1-2/N}$ in (\ref{k18}). In particular, at
the nanozone near the bottom of the potential, we have $N=2$,
$\hb'=1$, and (\ref{k22}) becomes 
\begin{equation}\label{k23}
\hb\bigg(n+\frac32\bigg)+\frac{\hb^{3/2}}{\sqrt{2}}
\bigg(\ol{z}\ol{\partial}^2-\frac12(\ol{z}^2-\ve_n)\ol{\partial}
+\frac{[{n}/{2}]}2\ol{z}\bigg)
+O(\hb^{2}).
\end{equation}

The model ordinary differential operator staying in (\ref{k23})
at the term $\hb^{3/2}$ determined the asymptotical properties
of the original Hamiltonian (\ref{k10}). By resolving the spectral
problem in the space $\cH_n$:
\begin{equation}\label{k24}
\frac1{\sqrt{2}}
\bigg(\ol{z}\ol{\partial}^2-\frac12(\ol{z}^2-\ve_n)\ol{\partial}
+\frac{[{n}/{2}]}2\ol{z}\bigg)
\varphi(\ol{z})
=\nu \varphi(\ol{z}),
\end{equation}
we obtain the eigenvalues $\nu=\nu_{n,k}$ and the eigenfunctions
(polynomials) $\varphi=\varphi_{n,k}$, where 
$k=0,\dots,[{n}/2]$.

\begin{theorem}%1
The asymptotics of the near bottom eigenvalues $\lambda$ and
the eigenfunctions $\psi$ of the Hamiltonian $\bH$ 
{\rm(\ref{k10})} is the following{\rm:}
\begin{align}
\lambda&=\hb(n+\frac32)+\hb^{3/2}\nu_{n,k}+O(\hb^2),
\label{k25}\\
\psi&=T(\varphi_{n,k})+O(\hb^{1/2}).
\nonumber
\end{align}
Here $n=0,1,\dots$ and $0\leq k\leq [{n}/2]$;
the numerating numbers $n,k$ are of order $O(1)$ as $\hb\to0$.
\end{theorem}

Equation (24), which gives the corrections $\nu_{n,k}$, 
and the coherent transform $T$ in (\ref{k25}) present 
the main difference of this resonance asymptotics from the
standard oscillatory approach. 

Of course, all the higher corrections of the asymptotics
(\ref{k25}) are calculated explicitly by simple perturbation
series.

If we go away from the nanozone near the bottom of the
potential (\ref{k10}) to some $N$th microzone, then we must
replace the model equation (\ref{k24}) by the equation 
\begin{equation}\label{k26}
\frac1{\sqrt{2}}
\bigg({\hb'}^2\ol{z}\ol{\partial}^2
-\frac{\hb'}2(\ol{z}^2-\ve_n\hb')\ol{\partial}
+\frac{\hb'[{n}/{2}]}2\ol{z}\bigg)\varphi
=\nu\varphi.
\end{equation}

Here $\hb'=\hb^{1-2/N}$ and $n\sim O(1/\hb')$. 
The number $N\geq3$, so $\hb'\to0$ and $n\gg1$. 
In this situation,
we can approximately solve Eq.~(\ref{k26}) using 
the technique of geometric coherent states over Lagrangian
submanifolds (classical trajectories) in the symplectic leaf
$\Omega$ developed in~\cite{kar24,kar25,kar10}. 

We now briefly describe the result.

In view of (\ref{k22}),
the classical Hamiltonian over $\Omega$ is given by the
coordinate function $\cA_3$ in the Poisson algebra (\ref{k12}).
Let us consider the energy levels of this Hamiltonian, that is
the closed curves $\Lambda\subset\Omega$ defined by 
\begin{equation}\label{k27}
\Lambda=\{\cA_3=\nu\}.
\end{equation}

Denote by $\Sigma$ a part of $\Omega$ bonded by the curve $\Lambda$
and consider the quantization condition
\begin{equation}\label{k28}
\frac1{2\pi\hb'}\int_{\Sigma}\omega_0=k+\frac12,\qquad k\in\bZ.
\end{equation}

This condition determines the discreet values $\nu=\nu_{n,k}$ in
(\ref{k27}). 
The numbers $\nu_{n,k}$ give the leading part of the eigenvalue
asymptotics in problem (\ref{k26}).
The asymptotics of the corresponding eigenfunctions $\varphi$ is
given by the integral over $\Lambda$:
\begin{equation*}\label{k29}
\varphi_{n,k}(\ol{z})=\frac1{\sqrt{2\pi\hb'}}
\int_{\Lambda}\Big(\sqrt{\dot z(t)}+O(\hb')\Big)
e^{\frac{i}{\hb'}\int^{t}_{0}\theta}\, K(\ol{z}z(t))\,dt.
\end{equation*}
Here $\{z=z(t)\}$ is the parametrization of points of $\Lambda$ 
by the time $t$ in the Hamiltonian system over $\Omega$ 
generated by $\cA_3$, the one-form $\theta$ 
is the primitive of the classical K\"ahlerian form~$\omega_0$, 
namely, $\theta=i\partial F_0$, 
and $K$ is the hypergeometric polynomial (\ref{k17}).

So, finally, we obtain the asymptotics of higher energy levels
of the original operator~$\bH$.

\begin{theorem}%2
Let $\hbar'=\hbar^{1-2/N}$, $N\geq3$.
The asymptotics of the eigenvalues of the operator
{\rm(\ref{k10})} with the numerating numbers
$n\sim1/\hb'$ is given by
\begin{equation}\label{k30}
\lambda=\hb\bigg(n+\frac32\bigg)+\hb^{3/N}\nu_{n,k}
+O(\hb^{4/N}).
\end{equation}
Here the values $\nu_{n,k}=\nu_{n,k}(\hbar')$ 
are determined by the quantization
condition {\rm(\ref{k28})}, where $0\leq k\leq[{n}/2]$. 
The asymptotics of the eigenfunctions corresponding to the
eigenvalues {\rm(\ref{k30})} is the following{\rm:}
\begin{equation}\label{k333}
\psi=\frac1{\sqrt{2\pi\hb'}}
\int_{\Lambda}\sqrt{\dot z(t)} \,
\exp\bigg\{\frac{i}{\hb'}
\int^{t}_{0}(\theta+O({\hb'}^{1/(N-2)}) )\bigg\}\, |z(t)\rangle\,dt,
\end{equation}
where $|z\rangle$ are coherent states of
the algebra {\rm(\ref{k8})} in the space $L^2(\bR^2)$ given by
{\rm(\ref{k19})}. 
All the corrections in the remainders 
in {\rm(\ref{k30})} and {\rm(\ref{k333})} 
are controlled by the higher-order terms in {\rm(\ref{k11})}
and are calculated explicitly.
\end{theorem}

By increasing the number $N=3,4,\dots$, 
we can go further and further 
away from the bottom point of the spectrum of $\bH$. 
But in any case the quantum number $n\sim\hb^{-(1-2/N)}$ in our
asymptotics can never reach the order $O(\hb^{-1})$. 
Such large numbers $n\sim \hb^{-1}$ correspond to
the energy levels of $\bH$ which are at the distance 
$O(1)$ from the bottom level. 
In this area, the behavior of the operator $\bH$ (\ref{k10}) 
is completely chaotic.  

The semiclassical parameter $\hbar'$ 
in (28) and (30) depends on the index~$N$  
of the microzone, and $\hbar'\gg \hbar$. 
One could take $\hbar' = \hbar$  
only in the chaos zone where the asymptotics (29) and (30) fail.

Note that arriving closer to the chaos zone, i.e., 
taking larger values of~$N$, one needs to take into account
higher-order terms in expression~(\ref{k9}) 
at least up to $f_{2N-2}$. 
All these terms give a contribution to the asymptotics of
eigenvalues (\ref{k30}) considered with the accuracy 
$o(\hb^2)$. This means that all the terms $V_j$ ($j=3,\dots,N$)
of the Taylor expansion (\ref{k2}) of the potential at 
the resonance bottom point contribute 
to the $o(\hb^2)$-asymptotics of
eigenvalues of the Schr\"odinger operator (\ref{k1}) 
in the $N$th microzone. 
But the geometry of the curve (\ref{k27}) and the leading 
quantization condition (\ref{k28}) are determined 
by the term~$V_3$ only.

\section{Resonance long-time evolution}%6

As a simple application of the above results, we describe the 
solution of the long-time evolution problem for the 
Schr\"odinger operator in the presence of resonance. 
Let us consider the Cauchy problem for the operator~(\ref{k10}):
\begin{align}
i\hb\frac{\partial \chi}{\partial t}
&=-\frac{\hb^2}{2}\Delta\chi
+\bigg(\frac12 x^2+2y^2+x^2y+\gamma x^4\bigg)\chi\,,
\label{k31}\\
\chi\bigg|_{t=0}
&=\chi^0({x}/{\sqrt{\hb}},{y}/{\sqrt{\hb}}),
\nonumber
\end{align}
where $\chi^0\in S(\bR^2)$. 
The initial data are localized in a nanozone $O(\sqrt{\hb})$
near the bottom point $x=y=0$ of the potential. Let us choose
the following time values:
\begin{equation}\label{k32}
t={\tau}/{\sqrt{\hb}},\qquad \tau\sim1.
\end{equation}
From the results described above, we obtain the asymptotics
\begin{equation}\label{k33}
\chi\approx\sum_{n\geq0}\exp\{-i(n+3/2)\tau/\sqrt{\hb}\}\,
\chi^{\tau}_{n}
({x}/{\sqrt{\hb}},{y}/{\sqrt{\hb}}),
\end{equation}
where $\chi^{\tau}_{n}$ is the solution of the evolution
equation in the $n$th irreducible representation of the algebra
(\ref{k8}): 
\begin{equation}\label{k34}
\bigg(-i\frac{\partial}{\partial \tau}
+\frac1{\sqrt{2}}\bA_3\bigg)\chi^{\tau}_{n}=0,\qquad 
\chi^{\tau}_{n}\bigg|_{\tau=0}=P_n(\chi^0).
\end{equation}    

Formula (\ref{k33}) demonstrates that the Schr\"odinger evolution
of the wave packet (\ref{k31}) in the long-time interval
(\ref{k32}) still keeps the packet to be localized.

The problem (\ref{k34}) after the coherent transform is reduced
to the evolution problem for the second order ordinary
differential operator 
on the left-hand side of (\ref{k24}).
This model problem does not have any small or large parameters 
and must be resolved exactly.

This is the {\it quantum nano-dynamics\/} describing the
Schr\"odinger evolution of wave packets localized at the
resonance bottom point.

For smaller time interval
\begin{equation}\label{k35}
t\sim O(\hb^{-1/N}),\qquad N\geq 3,
\end{equation}
we have, instead of (\ref{k34}), the problem with a small
parameter 
\begin{equation}\label{k36}
\bigg(-i\hb'\frac{\partial}{\partial \tau}
+\frac1{\sqrt{2}}\bA_3\bigg)\chi^{\tau}_{n}=0,
\qquad 
\hb'=\hb^{1-2/N},
\end{equation}
where $\bA_3$ is the $n$th irreducible representation is given by
the left-hand side of~(\ref{k26}). 
The evolution equation~(\ref{k36}) can be asymptotically solved 
by using the semiclassical approximation theory. So, 
in the time intervals like (\ref{k35}), the asymptotics of the
solution of the Cauchy problem (\ref{k31}) is computed
explicitly. 

Let us remark that the equations of the form 
(\ref{k36}) and more general equations which involve algebras of
the type (\ref{k6}) are related to some hidden 
geometry (classical and quantum). 
The global geometric analysis of these equations leads to
constructions of symplectic and quantum paths, to the symplectic
and quantum holonomy and curvature, to the translocation
operation, see~\cite{kar11,kar12,kar13}. 

It may be relevant to mention that 
in the absence of resonance, 
the Cauchy problem like (\ref{k31}) is asymptotically 
solved without any difficulties on the time interval 
$t\sim O(1/\hbar^\infty)$.

\section{Different types of resonance algebras}%7

It is known that there are many
different types of resonances. Their variety strongly
depends on stability or unstability of the first variation 
of the dynamical system generated by the leading 
part of the Hamiltonian. Resonances can be stable, 
unstable, neutral and combinations of these types.  
Respectively, the resonance algebra $M_0$ can be of compact,
noncompact, nilpotent type or be a mixture  
(direct product) of these types.  
We just briefly mention three simple examples.

\begin{example}[{\bf Landau model}] %Example 1.  
The well-known example of the neutral resonance is
given by the Hamiltonian of a charged particle moving
along a plane in a homogeneous magnetic field which is
perpendicular to the plane. 
In this case the first variation matrix of the Hamiltonian
field has zero eigenvalue which is twice degenerate. 
The resonance algebra~$M_0$ 
of functions commuting with the
Hamiltonian is nilpotent and just reduced to the three-dimensional
Heisenberg algebra.
\end{example}

\begin{example}[{\bf Inverted oscillator}] %Example 2.   
An example of the unstable resonance is given by the
inverted oscillator, i.e., the oscillator (2) which has not
a potential quadratic well, but a potential quadratic hill.
In  this case the oscillator frequencies $\alpha$, $\beta$
are imaginary, say,  $\alpha = \beta  = i$. 
Then the first variation matrix has two real eigenvalues
$\pm1$ which are twice degenerate.  The resonance algebra $M_0$
in this case is the enveloping of the Lie algebra $\su(1,1)$
which is of noncompact type.  

In the anisotropic version, where $\alpha = i$,
$\beta = 2i$, the resonance algebra~$M_0$ is given 
by quadratic commutation relations similar to~(8).
\end{example}

\begin{example}[{\bf Artificial magneto-atoms}] %Example 3.  
The resonances described in Section~3 are, of course,
stable. Their algebras are of compact type. 
One of physical models, where such resonances appear, 
is the so-called artificial atom  
(or $2$-dimensional quantum dot surrounded by electrons moving
in the plane). The resonance means the degeneracy of energy
levels of such an ``atom."  
This imply ``electron shells," ``filling numbers," etc.
and make such a system highly stable.

It is interesting to consider the artificial atom combined 
with a homogeneous magnetic field which is perpendicular
to the configuration plane.  
Let the potential well created by the central dot
be just quadratic with frequencies  
$\alpha =\beta = \omega_0$.  
Then we simply have the $2$-dimensional Fock model: 
oscillator plus magnetic field. 
Let us denote by $\omega_L$ half the Larmor frequency 
which is equal
to the magnetic field strength multiplied by the electron
charge and divided by the electron mass. One
can claim that the energy levels of the Fock 
model are degenerate if and only if
\begin{equation}
(\omega_L/\omega_0)^2=s^2/(k^2-s^2),
\label{8.1}
\end{equation}
where $s,\,k$  are some integers, $0\leq s < k$.  

Under this condition we obtain a {\it stable artificial
magneto-atom}.  
Its algebra of integrals of motion $M_0$ has four generators 
with polynomial commutation relations of type~(8),
and two Casimir functions. 
All components of the Poisson tensor (the right-hand side of the
commutation relations), except one, are linear. 
The nonlinear component is a polynomial of degree 
$d_0=l+m-1$, 
where~$l$, $m$ are coprime numbers such that
\begin{equation}
(k+s)/(k-s)=l/m.
\label{8.2}
\end{equation}

For different realizations of the magneto-atom,  that
is for different ratios (\ref{8.1}), 
one can have arbitrary numbers $l,\,m$  
in (\ref{8.2}) and make the degree~$d_0$  
of the algebra~$M_0$ be any number. 
For instance:

if $\omega_0 = 2\sqrt{2}\omega_L$,
then $d_0 = 2$,

if $\omega_0 = \sqrt{3}\omega_L$,  
then $d_0 = 3$,

if $\omega_0 = 2\sqrt{6}\omega_L$,
then $d_0 = 4$,

if $\omega_0 = \frac{4}{3}\omega_L$, 
then $d_0 = 4$,

if $\omega_0 = \frac{1}{2}\sqrt{5}\omega_L$, 
then $d_0 = 5$,  etc.

If there are higher terms $V_3,V_4,\dots$ in the potential
of the central dot, or a kind of an external field, 
then they determine a Hamiltonian over
the algebra~$M_0$ like (9). The spectrum and eigenfunctions 
of this Hamiltonian can be derived using the method described
in Sections~5 and~6. Thus one can calculate 
the splitting of the energy levels 
of this artificial magneto-atom 
with an arbitrary accuracy in~$\hbar$.

Note that if one replaced the central dot potential 
(the quadratic well) by the central Coulomb potential, 
i.e., considered the usual atom rather than
the artificial one, 
then its magnetic version might not exist at all.  
Indeed, it is well known that the model
``Coulomb potential plus magnetic field'' 
is chaotic;
there is no spectral degeneracy and therefore no
``magneto-atom."
\end{example}

In this example we obtain resonance algebras with
polynomial commutation relations of an arbitrarily high 
degree. The same growth of the degree is observed 
in {\it multi-dimensional\/} resonances. 
For three-dimensional oscillator 
with the resonance $1:2:3$
commutation relations in the algebra~$M_0$ 
are of the third degree, and so on.

\section{Nano-structures over resonance trajectories and tori}%8

Another class of problems which can be solved by the methods
described in Sections~5--7 is the asymptotics of spectral series
corresponding to closed stable trajectories or invariant
isotropic tori of the Hamiltonian dynamics in the presence of a
resonance between transversal Lyapunov frequencies or a
resonance  between elements of the Floquet holonomy group. 
The Maslov complex isotropic bundle theory \cite{Mas-Op,kar16}
in this situation has to be supplemented 
with an additional noncommutative algebraic and geometric
nano-structures like (\ref{k6}), (6a), (\ref{k8}), (\ref{k12}).

After application of the quantum
averaging method, we again reduce the problem to studying the
commutant $M_0$ of the leading part of the 
Floquet--Lyapunov transversal hamiltonian. 
The resonance means that the algebra $M_0$ is noncommutative, 
and we again can apply the quantum geometry technique to
calculate irreducible representations of $M_0$, coherent
states, etc. 
Here not only the hypergeometric type
but also the theta type \cite{kar9,kar19} quantum deformations 
can arise. 

In nano- and microzones near the trajectory (torus), 
in addition to the usual geometry~\cite{KarVor},
we obtain a
bundle with fiber-polynomial Poisson structures like 
(\ref{k12}), and the bundle of their symplectic leaves. 
The problem is reduced to model differential equations over
these bundles.  In the nanozone 
of order $O(\sqrt{\hb})$ these equations do not contain a small
parameter and must be resolved exactly, but in microzones 
of order $O(\hb^{1/N})$, $N\geq3$, 
the model equations can be effectively analyzed and 
explicitly solved using the global geometric semiclassical
technique. 

\begin{example}[{\bf Nano-electrodynamics\/}] %Example 4
This analysis can be applied, for instance, to the Helmholtz 
equation describing propagation of classical electromagnetic
waves in nonhomogeneous medium.  
Assume that the index of refraction of the medium has an
extremal, i.e., it takes the maximum value 
along a certain curve and does not degenerate  
in directions transversal to the curve. 
Then a neighborhood of this extremal becomes a wave channel. 
The Gaussian-type electromagnetic waves can propagate 
through this channel preserving their semiclassical localization
near the extremal (the channel axis).  

If the transversal frequencies are constant and are in resonance, 
then in nano- and microzones near the axis there appear certain
noncommutative algebraic (quantum) structures 
and Poisson structures 
of the type described  above. One can consider this 
phenomenon as a generation of a new type ``quantum"
polarization of the wave. The evolution of this polarization
is described by a model differential equation in the quantum
bundle over the extremal. Since the extremal is assumed to 
be stable, this polarization occurs to be of spin type. 
The ``spin" variables are running along symplectic leaves of 
algebras like (6a) and (8) which are diffeomorphic to a sphere.

Of course, these ``spin" variables have nothing to do with 
the actual physical spin. By the way, the situation is not
at the atomic scale, the scale  of the wave channel  can be
even macroscopic.

If the channel axis is not straight, but 
curvilinear, then the ``spin" degrees of freedom correlate
with its curvature and torsion. In this case one can speak
about a specific ``spin--orbit" interaction in the curvilinear
resonance wave channel.  
We detect existence of something like 
{\it quantum nano-electrodynamics\/} 
near resonance extremals of the index of refraction.
\end{example}

\section{Resonant quantization of wave mechanics}%9

One observes the following:

-- the resonance of Birkhoff or Lyapunov, or Floquet 
frequencies implies the appearance of specific
noncommutative algebras, in general, 
they are not Lie algebras;

-- such a noncommutative structure appears 
in equations of wave mechanics of various type,
for instance, in the Helmholtz and Maxwell equations, 
describing propagation of classical electromagnetic 
waves along resonance channels, and so the quantum behavior
can be inherent in classical wave systems too;

-- there is a nontrivial Poisson and symplectic geometry 
in nano- and micro- phase spaces near the resonance motion,
this resonance geometry is invisible in the usual classical
limit;

-- the structure of the phase space in nanozones near resonances
looses its classical behavior and occurs to be pure quantum, 
in particular, the resonance algebras have no classical Poisson 
limit in these zones, and so, the classical geometry must be
replaced by a quantum one;

-- to construct irreducible representations and coherent states
for the resonance algebras, one needs to introduce 
some new concepts to the quantization procedure, 
in particular, the quantum K\"ahlerian structures 
of hy\-per\-geo\-me\-tric- or theta-types;

-- in nanozones near the resonances, 
the given spectral problem is reduced to 
a model differential equation of order more than~$1$, 
in general;

-- the structure of the resonance algebra 
and of the accompanying quantum geometry are completely 
determined by the arithmetical proportion between
frequencies, one can consider these structures 
as a presentation of arithmetic in wave mechanics;

-- the micro- and nano- phase spaces related to the resonance
carry certain hidden dynamical geometry, classical and
quantum, in particular, quantum connections, quantum paths and
translocations; 

-- this dynamic geometry opens a new way to describe the
correlation and transition of modes in resonance clusters;

-- the general algebraic and geometric technique allows one 
to resolve the old problem of resonances in the
semiclassical approximation theory, in particular, in the
spectral and in the  long-time analysis.

On this way, we observe something interesting 
in spectral properties of the Schr\"odinger operator 
at resonance.  
For instance, in two dimensions, in the case of anisotropic
resonance, say $1:2$, the spectrum is given by a series 
in fractional powers $\hbar^{1/2}$ of the ``Planck constant''
(see (\ref{k25})). 

But, perhaps, the main find in the study of the
resonance problem is a class of very simple and fundamental
physical systems whose description involves 
{\it algebras with polynomial commutation relations\/} 
like~(\ref{k8}). This means that
algebras of such a non-Lie type are actually very common, 
although they were not clearly visible and
therefore their importance was not enough 
appreciated until now.  

We refer to~\cite{kar20}, and to the paper by Karasev \& Novikova 
in the given Collection, for another class of basic physical
examples where such types of algebras appear.

The ideas described above can imply a change 
in the viewpoint on the role of resonances in wave and quantum
problems.  First of all, 
systems with resonances possibly have to be considered not as
exceptional (like in the KAM theory), but as the basic ones
with additional perturbations which control deviations.
From the quantum attitude the following is clear:  
resonances imply high spectral multiplicities,  
this makes gaps between excited energy levels much wider,  
and therefore such systems  
become much more stable under perturbations. 
One can call this phenomenon a {\it resonant stabilization}. 
Note that all the fundamental atomic objects are resonant
stabilized. 
The world on the whole is supplied with 
enough amount of such a resonant stabilizers.

Secondly, the resonances can be a source of models with
noncommutative degrees of freedom (like spin or even with
non-Lie commutation relations). 
The notion of resonant stabilizer might be useful
for some extended physical concepts
like the string theory, as well as for models of quantum
computers and models used in modern micro- and nano-technologies. 

Among possible new applications there is the
opportunity to construct noncommutative optical channels
(see Example~4)
with an interference between ``spin" modes, 
with ``spin--orbit" effects, and with ``spin"  level 
tunneling  or crossing. 
Another application is the artificial magneto-atoms  
(see Example~3) with the opportunity 
to change their ``chemistry" just by varying resonance 
values of the magnetic field via formula~(37). 
The artificial atoms are also interesting from the viewpoint
of separation of their electron shells, which is closely
related to the phenomenon of quantization of the
configuration surface where they are located.

Next parts of the paper contain more detailed analysis of
resonance noncommutative algebras, as well the description of
several effects such as the resonance localization
and focusing, resonance traps, resonance adiabatic 
and ``spin" phenomena,  secondary resonances, etc.

\part{}

\setcounter{section}{-1}

\section{Introduction}%0

The resonance between frequencies of a Hamiltonian 
(or quantum) motion is a source of interesting structures and
effects in classical and wave mechanics \cite{kar14,kar15}.
It was observed in \cite{kar3,kar21} that resonances generate
noncommutative algebras which determine 
a specific ``quantum'' behavior 
of the system near a distinguished motion (or equilibria).
Here we investigate the resonance algebras  
in a more systematic way.

The first two sections below are purely algebraic.
We describe in detail the structure of the algebra 
of integrals of motion for the multidimensional harmonic
oscillator. 
If the frequencies of the oscillator are in a resonance, 
this algebra is noncommutative and commutation relations
are, in general, not of the Lie type.

In Section~3, we consider examples of resonance algebras 
for two-dimensional oscillators, 
including the inverted (unstable) oscillator, 
and the oscillator in a magnetic field. 

In Section~4, we describe the phenomenon of resonance 
precession in the case of two-dimensional oscillator. 
In Section~5, we discuss a general procedure of projecting 
the Poisson brackets and the appearance of the triple algebra
structure. We consider the examples of 1:1 and 1:2 resonances,
describe their triple algebras over the configuration plane, 
and demonstrate the way of solving the resonance precession
system in noncommuting space--time coordinates.

The Appendix contains a short list of formulas of the operator 
averaging method. 

\section{Commutant of the oscillator}%1

Let us consider a quadratic function $H_0$ 
on  $\bR^{2M}=\bR^M\times \bR^M$ and the corresponding
Hamiltonian vector field $J dH_0$, where 
$J=\left(\begin{matrix} 0&I\\-I&0\end{matrix}\right)$. 
Assume that the first variation matrix $JD^2H_0$ is stable, 
that is, it has only imaginary eigenvalues. 
Denote them by $\pm i\omega_l$ ($l=1,\dots,M$).
The function $H_0$, in suitable symplectic coordinates $q,\,p$,
can be written as 
\begin{equation}
H_0=\frac12\sum^{M}_{l=1}\omega_l(q^2_l+p^2_l).
\tag{1.1}
\end{equation}
One can call $H_0$ the {\it oscillator\/} Hamiltonian,
and call $\omega_l$ its {\it frequencies}.

We consider the algebra $F$ of all differential operators with
polynomial coefficients on $\bR^M$. An element of this algebra
is represented by the quantum oscillator
\begin{equation}
\wh{H}_0=\frac12\sum^{M}_{l=1}\omega_l
\Big(q^2_l-\hb^2\frac{\pa^2}{\pa q^2_l}-\hb\Big),
\tag{1.2}
\end{equation}
where $\hb>0$ is a fixed number.
Denote by $F_\omega$ the commutant of the oscillator (1.2) 
in the algebra~$F$.
By definition, 
$F_\omega$ consists of all elements from $F$ commuting
with $\wh{H}_0$. 
Obviously, $F_\omega$ is an algebra.

We also consider the Poisson algebra $\cF_\omega$ of functions
on $\bR^{2M}$ which are integrals of motion of the Hamiltonian
vector field $JdH_0$.  
Each element from $\cF_\omega$ has zero Poisson bracket 
with the function $H_0$ (1.1), i.e., Poisson commutes with
$H_0$. 
The Poisson algebra $\cF_\omega$ is the classical analog of the 
associative algebra $F_\omega$.

Our goal in this section is to study the properties of the
algebra $\cF_\omega$ (or $F_\omega$) related to the properties 
of the frequency system~$\omega$.  

We consider nonordered systems (sets) of real nonzero numbers
$\alpha=\{\alpha_1,\dots,\alpha_M\}$. 
We call $M$ the {\it length\/} of the system.

For any two systems $\alpha,\,\beta$, possibly, of different
lengths, there are naturally defined systems 
$\alpha\cup\beta$ and $\alpha\cap\beta$, as well as 
the notion of subsystem $\beta\subset\alpha$.

Also for any real number $c\neq0$, there are naturally defined
systems $c\cdot\alpha=\{c\alpha_1,\dots,c\alpha_M\}$.

A system $\alpha$ is said to be {\it positive\/} 
(nonnegative, or {\it integer\/}) 
if all elements $\alpha_l$ are positive  
(nonnegative, or integer).
We call $\alpha$ a {\it prime\/} system if all $\alpha_l$
are integer and coprime. 

A system $\alpha$ of length~$M$ is called {\it resonance\/}
if there is an integer system $m\neq0$ of length $M$ such that 
$\sum^{M}_{l=1}\alpha_l m_l=0$.

A system $\alpha$ is said to be {\it pure\/} 
if all its elements are mutually commensurable, i.e., 
if $\alpha_l/\alpha_k$ are rational for all $l,k$.

\begin{lemma}
Any pure system $\alpha$ can be written uniquely as 
$\alpha=\alpha_0\cdot n$, where $n$ is a prime system and 
$\alpha_0\in\bR$.
In particular, any pure system of length $>1$ is resonance.
\end{lemma}

The number $\alpha_0$ determined by Lemma~1.1 is called 
a {\it characteristic\/} of the pure system~$\alpha$.

\begin{lemma}
Let $\alpha,\,\beta$ be pure systems. 
Then the system $\alpha\cup\beta$ is pure if and only if 
the characteristics $\alpha_0$ and $\beta_0$ are commensurable. 
\end{lemma}

\begin{definition}
A pure subsystem of $\alpha$ which cannot be presented 
as a union of different pure subsystems will be called a {\it
component\/} of~$\alpha$. 
\end{definition}

\begin{lemma}
Different components do not intersect.
\end{lemma}

In summary, this results in the following statement.

\begin{proposition}
Each system of nonzero numbers can be uniquely resolved 
into a disjoint union of components. The characteristics of the
components are incommensurable. The system is resonance if and
only if at least one of its components has length greater than~$1$.
\end{proposition}

Now we pass from systems of numbers to algebras.

Let $B,\,C$ be two subalgebras in an algebra $F$. 
The {\it algebraic envelope\/} $B\dot{+}C$ is a subalgebra
in~$F$ generated by all sums and products of elements from~$B$
and~$C$.  
The {\it commutator\/} $[B,C]$ is a subalgebra in~$F$ generated
by commutators $[b,c]$ of arbitrary elements $b\in B$ and $c\in C$.
The subalgebras $B$ and $C$ commute if $[B,C]=0$.

A subalgebra $B$ is said to be {\it pure\/} if it cannot be
presented as an algebraic envelope of different commuting
subalgebras. 

All pure commutative subalgebras are one-dimensional (i.e.,
generated by a single, non-unity, element).

\begin{definition} 
A pure subalgebra in~$F$, 
which is not contained in any larger pure subalgebra, 
will be called a {\it component\/} of the algebra~$F$.
\end{definition}

Now let us return to the commutant $F_\omega$ of the oscillator
$H_0$ (1.2). The previous definitions and statements allow us 
to formulate the following result.

\begin{theorem}
The commutant $F_\omega$ of the oscillator is uniquely
presented as an algebraic envelope of its commuting components
\begin{equation}
F_\omega=F^{(1)}\dot{+}\dots\dot{+}F^{(L)},
\qquad
[F^{(j)},F^{(s)}]=0.
\tag{1.3}
\end{equation}
This resolution corresponds to the resolution of the frequency
system $\omega$ into its components 
{\rm(}see Proposition~{\rm1.1):}
\begin{equation}
\omega=\bigcup^{L}_{j=1}\omega^{(j)},\qquad 
\omega^{(j)}\cap\omega^{(s)}=\emptyset.
\tag{1.4}
\end{equation}
The characteristics $\omega^{(j)}_0$ of these components are
incommensurable. 

The algebra $F_\omega$ is noncommutative if and only if the
frequency system $\omega$ is resonance. This happens if and only
if at least one of the components $\omega^{(j)}$ has length 
greater than $1$ {\rm(}and so, the corresponding component
$F^{(j)}$ is greater than one-dimensional\/{\rm)}. 
\end{theorem}

\begin{proof}
Note that the Hamiltonian (1.1) can be written in the form 
\begin{equation}
H_0=\sum^{M}_{l=1}\omega_l\oz_l z_l,\qquad
z_l=\frac1{\sqrt{2}}(q_l+ip_l).
\tag{1.1a}
\end{equation}
The oscillator (1.2) is respectively written as 
\begin{equation}
\wh{H}_0=\sum^{M}_{l=1}\omega_l\hat{z}^*_l \hat{z}_l,\qquad
\hat{z}_l=\frac1{\sqrt{2}}\Big(q_l+\hb \frac{\pa}{\pa q_l}\Big),
\tag{1.2a}
\end{equation}
where $\hat{z}^*$ and $\hat{z}$ are ordered in the Wick sense
(all $\hat{z}_l$ on the right, all $\hat{z}^*_l$ on the left),
the asterisk means the adjoint operator in $L^2(\bR^M)$.

For any Wick ordered polynomial $g$ in the operators $\hat{z}^*$,
$\hat{z}$ we have 
\begin{align*}
[\wh{H}_0,g(\hat{z}^*,\hat{z})]
&=-i\hb\{H_0,g\}(\hat{z}^*,\hat{z})
=\hb(\pa H_0 \opa g-\opa H_0\pa g)(\hat{z}^*,\hat{z})\\
&=\hb \sum^{M}_{l=1}\omega_l
(\oz_l\opa_lg-z_l\pa_lg)(\hat{z}^*,\hat{z}).
\end{align*}
Therefore (see \cite{Gust}) the commutant $F_\omega$ is generated
by the elements 
\begin{equation}
\hat{g}_k=(\hat{z}^*)^{k_-}\hat{z}^{k_+},\qquad
\omega\circ k_+=\omega\circ k_-.
\tag{1.5}
\end{equation}
Here the systems $k_+,k_-$ are integer nonnegative and have the
same length $M$ as the frequency system $\omega$, 
the notation~$\circ$ is used for the scalar product in $\bZ^M_+$:
\begin{equation}
\alpha\circ \beta=\sum^{M}_{l=1}\alpha_l\beta_l.
\tag{1.6}
\end{equation}

Let the frequency system $\omega$ be resolved by its components
(1.4). By Lemma~1.1, each component can be written as 
$\omega^{(j)}=\omega^{(j)}_0\cdot n^{(j)}$, where
$n^{(j)}$ is an integer system. Then scalar products like (1.6)
are represented as the sums
$$
\omega\circ\alpha=\sum^{L}_{j=1}
\omega^{(j)}_0 n^{(j)}\circ\alpha^{(j)},
$$
where $\alpha^{(j)}$ are disjoint subsystems of the system
$\alpha$, that is, $\alpha=\bigcup^L_{j=1}\alpha^{(j)}$,
$\alpha^{(j)}\cap \alpha^{(s)}=\emptyset$. 
Thus the basic equation in (1.5) is transformed into
$$
\sum^{L}_{j=1}\omega^{(j)}_0 n^{(j)}\circ(k^{(j)}_+-k^{(j)}_-)=0.
$$
Since all the characteristics $\omega^{(j)}_0$ are mutually 
incommensurable, this equation is equivalent to 
\begin{equation}
n^{(j)}\circ(k^{(j)}_+ - k^{(j)}_-)=0,\qquad
\forall j=1,\dots,L.
\tag{1.7}
\end{equation}
This means that the $j$th component $F^{(j)}$ of the commutator
$F_\omega$ is generated by the operators $\hat{g}_{k^{(j)}}$
of type (1.5),
where the integers systems $k^{(j)}_+$, $k^{(j)}_-$ 
obey the $j$th equation in (1.7).

From Proposition 1.1, the frequency system is resonance 
if and only if for some $j$ the system $n^{(j)}$ has length
greater than~$1$. In this case, 
there are at least two solutions of the $j$th equation in (1.7) 
and the corresponding monomials $g_{k^{(j)}}$ and $g_{k^{\prime(j)}}$ 
whose mutual Poisson bracket is not zero, 
and therefore the algebra $F^{(j)}$ is noncommutative.
\end{proof}

\begin{remark}
The classical analog of Theorem 1.1 is obvious:
just replace the word ``algebra'' by the words ``Poisson
algebra'' and replace the commutator in (1.3) by the Poisson
brackets. 
\end{remark}

Theorem 1.1 means that in the nonresonance case, where all
frequencies $\omega_l$ are incommensurable, the commutant
$F_\omega$ is trivial and generated by the commuting elements 
\begin{equation}
\hat{g}_{I_l}\od \hat{z}^*_l\hat{z}_l,\qquad l=1,\dots,M.
\tag{1.8}
\end{equation}
In the resonance case the study of the algebra $F_\omega$ is
reduced to the study of its nontrivial components.

Therefore, in what follows, we consider only the one-component
situation, that is, we assume that the frequency system has the
form 
$$
\omega=\{n_1,\dots,n_M\},\qquad M\geq2,
$$
where $n_l$ are integers. Moreover, without loss of generality, 
we can assume that $n_l$ are positive, since all minus signs 
can be included into the systems of numbers $k_+,\,k_-$ in (1.7).

Let us consider the oscillator with positive integer frequencies:
\begin{equation}
\wh{H}_0=\sum^{M}_{l=1} n_l\hat{z}^*_l\hat{z}_l,\qquad
n_l\in\bN.
\tag{1.9}
\end{equation}
Its commutant $F_n$ will be called a {\it resonance algebra}. 
This algebra is generated by the operators $\hat{g}_k$ (1.5),
where 
\begin{equation}
n\circ (k_+ - k_-)=0,\qquad k_+,k_-\in\bZ^M_+.
\tag{1.10}
\end{equation}

Let us denote by $R_n$ the resonance set of all solutions 
$k=(k_+,k_-)$ of Eq.~(1.10). 
Evidently, $R_n$ is an Abelian subsemigroup in $\bZ^M_+\times
\bZ^M_+$.  We can consider the secondary subsemigroup
\begin{equation}
R^{(1)}_n=R_n+R_n,
\tag{1.11}
\end{equation}
and take its complement
\begin{equation}
\cR_n\od R_n\setminus R^{(1)}_n.
\tag{1.12}
\end{equation}

This later set consists of those solutions of Eq.~(1.10) which
cannot be presented as a sum of two other solutions.
The simplest example of an element from $\cR_n$ is
\begin{equation}
I_l\od\big(\,
\underbrace{0,\dots,0}_{l-1},1,\underbrace{0,\dots,0}_{M-l};
\underbrace{0,\dots,0}_{l-1},1,\underbrace{0,\dots,0}_{M-l}\,\big).
\tag{1.13}
\end{equation}

Such $I_l\in \cR_n$ we call {\it primitive\/} resonance elements.
A generic element from $\cR_n$ will be called a {\it minimal\/}
resonance element.  

\begin{theorem} 
The number of minimal resonance elements is finite. Each
solution of the resonance equation {\rm(1.10)} can be
represented as a linear combination of minimal resonance
elements with integer positive coefficients.
\end{theorem}

\begin{proof}
Let $k=(k_+,k_-)$ be a solution of (1.10) and let it be minimal,
i.e., $k\in\cR_n$. We claim that for any $l=1,\dots,M$
\begin{equation}
k_{+l}\leq \sum n,\qquad k_{-l}\leq \sum n.
\tag{1.14}
\end{equation}
Here the sum $\sum$ is taken over all indices $1,\dots,M$.
To prove (1.14), we, on the contrary, assume that (1.14) 
does not hold, say,
\begin{equation}
k_{+1}>\sum n.
\tag{1.15}
\end{equation}
Then we have $k_{+1}>\max n$. It follows from this inequality 
that all $k_{-l}$ do not exceed $n_1$. Indeed, if, on the contrary,
there exists some $k_{-l}$ such that 
\begin{equation}
k_{-l}>n_1,
\tag{1.16}
\end{equation}
then we obtain two different solutions $k'$ and $k''$ of Eq.~(1.10),
where 
$$ 
k'_+=(k_{+1}-n_l, k_{+2},\dots,k_{+M}),\qquad 
k'_-=(k_{-1}, k_{-l}-n_1,\dots,k_{-M}),
$$ 
and 
$$ 
k''_+=(n_l,0,\dots,0),\qquad
k''_-=(0,\dots,n_1,\dots,0).
$$
Their sum equals the original solution: $k'+k''=k$.

This contradicts the assumption that $k\in\cR_n$. Thus (1.16)
fails, and we conclude that 
$$
k_{-l}\leq n_1,\qquad\forall l=1,\dots,M.
$$
Then we derive 
$$
k_{+1}n_1\leq k_+\circ n=k_-\circ n
\leq \max_l (k_{-l})\sum n\leq n_1\sum n,
$$
and so $k_{+1}\leq \sum n$. This contradicts the assumption
(1.15), and hence we have proved (1.14).

It follows from (1.14) that the number of elements in $\cR_n$
does not exceed $2M\sum n$.

The second statement of the theorem is evident and related to
the fact that the multiplicative semigroup $\bN$ acts on $R_n$:
$$
\mu\cdot(k_+,k_-)\od (\mu k_+,\mu k_-),\qquad \mu=1,2,\dots,
$$
and the additive semigroup $\bZ^M_+$ acts too:
$$
m\cdot (k_+,k_-)\od (k_+ +m, k_- +m),\qquad m\in \bZ^M_+.
$$
The proof of the theorem is complete.
\end{proof}

\section{Polynomial relations in resonance algebras}%2

The resonance Poisson algebra $\cF_n$ is generated by the
functions $g_k=z^{k_+}\oz^{k_-}$ on $\bR^{2M}$, 
where systems of nonnegative integers 
$(k_+,k_-)\in \bZ^M_+\times\bZ^M_+$
obey Eq.~(1.10). 
It is easy to derive the Poisson bracket between a pair of
generators 
\begin{equation}
\{g_k,g_r\}_{\bR^{2M}}=i\sum^{M}_{l=1}[k,r]_l\, g_{k+r-I_l}.
\tag{2.1}
\end{equation}
Here $I_l$ is the primitive element (1.13), 
the bracket $[k,r]\in\bZ^M$ is defined by 
\begin{equation}
[k,r]=k_+r_- - k_-r_+,
\tag{2.2}
\end{equation}
where the multiplication in $\bZ^M_+$ is naturally determined by
its multiplicative semigroup structure.

One can now consider an abstract vector space  
with the set of complex coordinates $G_k$
such that $\ol{G}_k = G_{k^*}$.
Here $k\in R_n$ and the involution $k\to k^*$ is defined by 
\begin{equation}
k^*_+=k_-,\qquad k^*_-=k_+.
\tag{2.3}
\end{equation}
Let us set the bracket between a pair of these coordinates 
just by mimicing the bracket (2.1):
\begin{equation}
\{G_k,G_r\}\od i\sum^{M}_{l=1}[k,r]_l\, G_{k+r-I_l},\qquad 
k,r\in R_n.
\tag{2.4}
\end{equation}

\begin{proposition}
The bracket {\rm(2.4)} is skew-symmetric and obeys the Jacobi
identity, i.e., it is a Poisson bracket. This Poisson structure
is consistent with the involution 
\begin{equation}
\ol{\{G_k,G_r\}}=\{\ol{G}_k,\ol{G}_r\}.
\tag{2.5}
\end{equation}
\end{proposition}

\begin{proof}
The skew-symmetry is evident. To check the Jacobi identity, let
us compute the double bracket
$$
\{\{G_k,G_r\},G_s\}
=\sum^{M}_{l,j=1}\,[k,r]_l [k+r-I_l,s]_j\,G_{k+r+s-I_l-I_j}.
$$
On the left and on the right of this relation, one has to take
the cyclic sum $\fS$ over permutations of the indices $k,r,s$
and to prove that the right-hand side is identically zero.

On the right, we obtain cyclic sums of two types 
$$
\sum_{l,j}\big(\fS[k,r]_l[k,s]_j+\fS[k,r]_l[r,s]_j\big)
G_{k+r+s-I_l-I_j}
$$
and 
$$
\sum_l\big(\fS[k,r]_l[s,I_l]_l\big)G_{k+r+s-2I_l}.
$$

By the direct computation one can verify that 
the numerical coefficients in round brackets in both of these
expressions are identically zero (in the first expression one
even need not open the brackets $[\,,\,]$, 
but in the second expression the brackets $[\,,\,]$ 
must be opened using (2.2) to see that all summands cancel each
other). 
\end{proof}

The linear Poisson tensor on the right-hand side 
of relations (2.4) determines a Lie--Poisson structure on the
resonance algebra $\cF_n$.
But this structure is infinite dimensional. Our goal now is to
replace it by a finite dimensional (finitely generated) Poisson
algebra structure.

First of all, we note that the original brackets (2.1) can be
written in the following form
\begin{equation}
\{g_k,g_r\}_{\bR^{2M}}=i g_k g_r \Phi_{kr}(g_I).
\tag{2.6}
\end{equation}
Here the generators $g_I$ correspond to the primitive
elements~$I$ (1.13)
\begin{equation}
g_{I_l}=|z_l|^2,\qquad l=1,\dots,M,
\tag{2.7}
\end{equation}
and the functions $\Phi_{kr}$ are defined by
\begin{equation}
\Phi_{kr}(\lambda)\od \sum^{M}_{l=1} 
\frac{[k,r]_l}{\lambda_l},\qquad 
\lambda=(\lambda_1,\dots,\lambda_M).
\tag{2.8}
\end{equation}

Thus, we see from (2.6) that the bracket of two generators of
the algebra $\cF_n$ is proportional to the product of these
generators with a coefficient which is a combination of the
inverse primitive generators (2.7).
Although this coefficient has singularities, the whole
right-hand side of (2.6) is, of course, smooth.

We consider a domain where all 
$z_l\neq0$ in order to avoid singularities. Anyway, the final
result will not contain any singularities.

Let us consider an abstract set of complex coordinate functions
$\cA_k$, where $k\in R_n$, which is consistent with
the involution (2.3): $\ol{\cA}_k=\cA_{k^*}$. 
The coordinate $\cA_I$ is said to be {\it primitive\/} if the
element $I\in R_n$ is primitive, and $\cA_k$ is said to be 
{\it minimal\/} if the element $k\in \cR_n$ is minimal.

Let us introduce a bracket between the coordinates $A_k$ by
mimicing the structure of the bracket (2.6):
\begin{equation}
\{\cA_k,\cA_r\}\od i\cA_k\cA_r\Phi_{kr}(\cA_I).
\tag{2.9}
\end{equation}
It follows from (2.9) that the primitive coordinates $\cA_I$ are
in involution with each other:
\begin{equation}
\{\cA_{I_l},\cA_{I_j}\}=0,\qquad \forall l,j=1,\dots,M,
\tag{2.10}
\end{equation}
since $[I_l,I_j]=0$ (see definitions (2.2) and (1.13)).

From (2.9) we also obtain  
\begin{equation}
\{\cA_k,\cA_{I_j}\}=i(k_+ - k_-)_j \cA_k,
\tag{2.11}
\end{equation}
since $[k,I_j]_l=(k_+ - k_-)_j \delta_{jl}$.

Some other useful properties of the operation (2.2) 
we collect in the following statement.

\begin{lemma} 
The operation {\rm(2.2)} 
$\bZ^M_+\times\bZ^M_+\to\bZ^M$ has the properties
\begin{align}
&\text{\rm(a)}\quad [k+s,r]=[k,r]+[s,r],
\nn\\
&\text{\rm(b)}\quad [\mu\cdot k,r]=\mu\cdot [k,r],\qquad \mu\in \bZ_+,
\nn\\
&\text{\rm(c)}\quad [k,r]=[r^*,k^*],
\tag{2.12}\\
&\text{\rm(d)}\quad [k,r]=-[r,k],
\nn\\
&\text{\rm(e)}\quad \fS_{k,r,s}[k,r](s_+ - s_-)=0.
\nn
\end{align}
\end{lemma}

\begin{proof}
The only non-evident relation is the last one. 
Computing the cyclic sum, one obtains
\begin{align*}
\fS_{k,r,s}[k,r](s_+ - s_-) 
&=(k_+r_- - k_-r_+)(s_+ - s_-)
+(r_+s_- - r_-s_+)(k_+ - k_-)\\
&\qquad +(s_+k_- - s_-k_+)(r_+ - r_-).
\end{align*}
This is identical to zero.
\end{proof}

\begin{corollary}
The bracket defined by {\rm(2.9)} is a Poisson bracket. 
It is consistent with the involution
\begin{equation}
\ol{\{\cA_k,\cA_r\}}=\{\ol{\cA_k},\ol{\cA_r}\}.
\tag{2.13}
\end{equation}
\end{corollary}

\begin{proof}
The skew-symmetry $\{\cA_k,\cA_r\}=-\{\cA_r,\cA_k\}$ follows from
(2.12)\,(d).To check the Jacobi identity we derive 
\begin{align*}
\{\{\cA_k,\cA_r\},\cA_s\}
&=i\{\cA_k,\cA_s\}\cA_r \Phi_{kr}(\cA_I)
+i\{\cA_r,\cA_s\}\cA_k \Phi_{kr}(\cA_I)\\
&\qquad
+i\cA_k\cA_r\{\Phi_{kr}(\cA_I),\cA_s\}.
\end{align*}
Applying once more (2.9), as well (2.11), we obtain 
$$
\{\{\bA_k,\bA_r\},\bA_s\}
=-\cA_k \cA_r \cA_s\big[(\Phi_{ks}+\Phi_{rs})\Phi_{kr}
+\Gamma_{krs}\big],
$$
where
$$
\Gamma_{krs}\od\sum^{M}_{l=1}
\frac{[k,r]_l(s_+ - s_-)_l} {\cA_{I_l}^2}.
$$ 
It follows from (2.12)\,(e) that 
$\fS_{k,r,s}\Gamma_{krs}=0$.
From the skew-symmetry of $\Phi_{kr}$ we also obtain 
$\fS_{k,r,s}(\Phi_{ks}+\Phi_{rs})\Phi_{kr}=0$. 
Relation (2.13) follows from (2.12)\,(c).
\end{proof}

Now we can proof the key lemma.

\begin{lemma} 
If an element $k\in R_n$ obeys the expansion
\begin{equation}
\sum_{r\in\cR_n}m_r\cdot r=k,\qquad m_r\in\bZ_+,
\tag{2.14}
\end{equation}
then
\begin{equation}
\prod_{r\in\cR_n}\cA_r^{m_r}=C \cA_k,
\tag{2.15}
\end{equation}
where $C=C(\cA)$ is a Casimir function for the Poisson bracket 
{\rm(2.9)}. 
\end{lemma}

\begin{proof} 
We derive
$$ 
\Big\{\prod_r \cA_r^{m_r},\cA_s \Big\}
=\sum_{r} m_r \cA^{m_r-1}_r \{\cA_r,\cA_s\}
\prod_{r'\neq r}\cA_{r'}^{m_{r'}}.
$$
By taking the bracket $\{A_r,A_s\}$ from the definition (2.9)
and by using (2.8), (2.12\,a), and (2.14), 
we transform the previous formula as follows: 
\begin{align*}
\Big\{\prod_r \cA_r^{m_r},\cA_s \Big\}
&=i\cA_s\Big(\sum_{r}m_r\Phi_{rs}\Big)
\prod_{r'} \cA_{r'}^{m_{r'}}
\\
&=i\cA_s\Big(\sum^{M}_{l=1}\frac1{\lambda_l}
 \Big[\sum_r m_r\cdot r,s\Big]\Big)\prod_{r'}\cA_{r'}^{m_{r'}}
\\
&=i\cA_s\Phi_{ks}\prod_{r}\cA_{r}^{m_{r}}
=\frac{\{\cA_k,\cA_s\}}{\cA_k}\prod_{r}\cA_{r}^{m_{r}}.
\end{align*}
Therefore,
$$
\Big\{ \frac1{\cA_k}\prod_{r}\cA_{r}^{m_{r}},\cA_s\Big\}=0,
\qquad \forall s\in R_n,
$$
and hence the function $C=\frac1{\cA_k}\prod_{r}\cA_{r}^{m_{r}}$ 
is a Casimir function.
\end{proof}

Using this lemma, one can transform the right-hand side of (2.9)
into a form which does not have a singularity.

From now on, we consider only minimal generators $\cA_k$, 
$k\in \cR_n$.
The set of minimal elements $\cR_n$ can be divided into the
disjoint union
\begin{equation}
\cR_n=\cR^0_n\cup \cR^\#_n,
\tag{2.16}
\end{equation}
where $\cR^0_n$ consists of all primitive elements (1.13) 
and $\cR^\#$ consists of {\it clean\/} elements $k$ which obey
\begin{equation}
k_+ k_-=0.
\tag{2.17}
\end{equation}

Let us fix $l\in(1,\dots,M)$. For each pair of clean elements 
$k,r$ such that the numbers $k_{+l}$, $k_{-l}$, $r_{+l}$, $r_{-l}$ 
are not all zero, there is an expansion of the element $k+r-I_l$ 
via the minimal ones 
\begin{equation}
k+r-I_l=\sum_{s\in \cR_n} m^{l}_{k,r,s}\cdot s,\qquad
m^{l}_{k,r,s}\in \bZ_+.
\tag{2.18}
\end{equation}
Here we have $m^{l}_{k,r,s}=m^{l}_{r,k,s}=m^{l}_{k^*,r^*,s^*}$.

From Lemma~2.2 we conclude that there are some Casimir
functions $C^l_{k,r}$ such that 
\begin{equation}
\cA_k \cA_r=C^l_{k,r} \cA_{I_l}\prod_{s\in \cR_n}\cA_s^{m^l_{k,r,s}}.
\tag{2.19}
\end{equation}
Here we have $\ol{C^l_{k,r}}=C^l_{k^*,r^*}$ and 
$C^l_{k,r}=C^l_{r,k}$.

Then the Poisson bracket (2.9) reads
\begin{equation}
\{\cA_k,\cA_r\}=i\sum^{M}_{l=1}C^l_{k,r}[k,r]_l
\prod_{s\in\cR_n}\cA_s^{m^{l}_{k,r,s}}
\tag{2.20}
\end{equation}
for any two clean elements $k,r$. 

Formulas (2.20), together with (2.10), (2.11) 
provide a nonsingular expression for the Poisson
brackets of two minimal coordinates via the minimal coordinates.
From Theorem~1.2 we know that the number of minimal coordinates
$N=\#(\cR_n)$ is finite. 

In view of the involution property 
$\ol{\cA}_k=\cA_{k^*}$, 
the real $\Re(\cA_k)$ and imaginary $\Im(\cA_k)$ parts of minimal
coordinates belong to the same Poisson algebra.
Formula (2.13) demonstrates that the finite dimensional Poisson 
tensor which one sees in (2.20), (2.10), and (2,11) is actually
a real tensor on the constraint surface determined by Eqs.~(2.19).

Thus we obtain a finite dimensional Poisson algebra generated by
the minimal $\cA_k$ ($k\in\cR_n$).

In the realization of this algebra by functions on the phase
space $\bR^{2M}$,
\begin{equation}
\cA_k\to g_k=z^{k_+}\oz^{k_-},
\tag{2.21}
\end{equation}
all the Casimir elements $C$ in (2.15) are, of course, equal
to~$1$.  
So, from (2.19) we obtain a distinguished constraint surface
given by the equations
\begin{equation}
\cA_k\cA_r=\cA_{I_l}\prod_{s\in\cR_n}\cA_s^{m^l_{k,r,s}},\qquad 
k,r\in\cR^\#_n.
\tag{2.22}
\end{equation}
On this particular surface we observe the following Poisson
structure 
\begin{equation}
\{\cA_k,\cA_r\}=i\sum^{M}_{l=1}[k,r]_l
\prod_{s\in\cR_n} \cA_s^{m^l_{k,r,s}},
\qquad 
k,r\in\cR^\#_n.
\tag{2.23}
\end{equation}

Note that by a small change in the notation 
one can make the bracket $[\,,\,]$ in (2.23) 
be the actual bracket on the set of integers $\bZ^M$. 
To do this, we map $\bZ^M_+\times \bZ^M_+$ into $\bZ^M$:
$$
k\to k_+ - k_-.
$$
This map is one-to-one on the subset of clean elements obeying
(2.17). The inverse map is given by 
$$
\alpha\to(\alpha_+,\alpha_-),
$$
where $\alpha_{+l}$ and $\alpha_{-l}$ are the positive and
negative parts of the number $\alpha_l$. 
The resonance equation (1.10) reads 
\begin{equation}
n\circ \alpha=0.
\tag{2.24}
\end{equation}
The condition (2.17) is transformed to 
\begin{equation}
\alpha_+\alpha_-=0.
\tag{2.17a}
\end{equation}
The bracket (2.2) and the involution (2.3) are naturally
transported to $\bZ^M$ as
\begin{equation}
[\alpha,\beta]=\alpha_+\beta_- - \alpha_-\beta_+,\qquad
\alpha^*=-\alpha.
\tag{2.25}
\end{equation}
Properties (2.12)\,(b,c,d) hold for this bracket as well.

Now the analog of the property (2.12)\,(e) can be written as 
\begin{equation}
\fS_{\alpha,\beta,\gamma}[\alpha,\beta]\gamma=0.
\tag{2.26}
\end{equation}
The analog of property (2.12)\,(a) reads
\begin{equation}
[\alpha+\beta,\gamma]=[\alpha,\gamma]+[\beta,\gamma]
+(\alpha\cplus\beta)\gamma,
\tag{2.27}
\end{equation}
where 
\begin{equation}
\alpha\cplus\beta\od 
\min(\alpha_+ + \beta_+,\,\alpha_- + \beta_-)\in\bZ^M_+.
\tag{2.28}
\end{equation}

Thus the Jacobi-like property (2.26) looks natural, but the
linearity-like property (2.27) contains an ``anomaly.''
This anomaly itself obeys the following properties:
\begin{gather*}
\alpha\cplus\beta
=\alpha_+ + \beta_+ - (\alpha+\beta)_+
=\alpha_- + \beta_- - (\alpha+\beta)_-,
\\
(\alpha\cplus\beta)+((\alpha+\beta)\cplus\gamma)
=(\alpha\cplus(\beta+\gamma)) +(\beta\cplus\gamma).
\end{gather*}

The expansion (2.18) is replaced by 
$$
\alpha+\beta
=\sum_{\gamma\in\Gamma_n}\mu_{\alpha,\beta,\gamma}\cdot\gamma
+(\alpha\cplus\beta),\qquad  \alpha,\beta\in\Gamma_n.
$$
Here $\mu_{\alpha,\beta,\gamma}\in\bZ_+$, 
and 
by $\Gamma_n\subset\bZ^M$ we denote the set of solutions 
of the resonance equation (2.24).

The constraint equations (2.22) are replaced by 
\begin{equation}
\cA_\alpha \cA_\beta=\cA^{\alpha\cplus\beta}_I
\prod_{\gamma\in\Gamma_n}\cA_\gamma^{\mu_{\alpha,\beta,\gamma}},
\tag{2.29}
\end{equation}
and the involution condition is written as 
$$
\ol{\cA}_\alpha=\cA_{-\alpha}\quad (\alpha\in \Gamma_n),
\qquad  \ol{\cA}_I=\cA_I.
$$
Now the Poisson brackets are finally read
\begin{equation}
\begin{aligned}
\{\cA_\alpha,\cA_\beta\}&=i\cPhi_{\alpha\beta}(\cA_I)
\prod_{\gamma\in\Gamma_n} \cA_\gamma^{\mu_{\alpha,\beta,\gamma}},
\qquad \alpha,\beta\in\Gamma_n,
\\
\{\cA_\alpha,\cA_{I_j}\}&=i\alpha_j \cA_\alpha,\qquad
\{\cA_{I_l},\cA_{I_j}\}=0,\qquad
l,j=1,\dots,M.
\end{aligned}
\tag{2.30}
\end{equation}
Here the functions $\cPhi_{\alpha\beta}$ are defined by 
$$
\cPhi_{\alpha\beta}(\lambda)
=\lambda^{\alpha\cplus\beta}\sum^{M}_{l=1}
\frac{[\alpha,\beta]_l}{\lambda_l},\qquad
\lambda\in\bR^M,
$$
where the bracket $[\alpha,\beta]$ and the operation 
$\alpha\cplus\beta$ are given by (2.25) and (2.27). 

Note that for $\beta=-\alpha$, the constraint (2.29) and the
bracket (2.30) are especially simple:
\begin{align}
|\cA_\alpha|^2&=\cA^{\alpha_+ + \alpha_-}_I,
\tag{2.29a}\\
\{\cA_\alpha,\ol{\cA}_\alpha\}&=i \cA^{\alpha_+ + \alpha_-}_I
\sum^{M}_{l=1}\frac{\alpha^2_+ - \alpha^2_-}{\cA_{I_l}}.
\tag{2.30a}
\end{align}
Also note that 
\begin{equation}
[\alpha,\beta]=0\qquad\text{if and only if}\quad \alpha\beta=0\quad
(\text{where}\quad \alpha,\beta\in\Gamma_n),
\tag{2.31}
\end{equation}
and in this case $\{\cA_\alpha,\cA_\beta\}=0$.

\begin{theorem}
The algebra $\cF_n$ of functions on $\bR^{2M}$ which are in
involution with the oscillator 
$H_0=\frac12\sum^{M}_{l=1}n_l(q^2_l+p^2_l)$, $n_l\in\bN$, 
is a finite generated involutive Poisson algebra. 
After the identification 
$\cA_\alpha\leftrightarrow z^{\alpha_+}\ol{z}^{\alpha_-}$, 
$\cA_{I_j}\leftrightarrow |z|^{2I_j}$, 
the constraints and Poisson brackets 
of this algebra are of polynomial type and 
given by {\rm(2.29), (2.30)}. 
These brackets are arithmetic: all the structural
coefficients are integer numbers. The functional dimension $N$
of the algebra $\cF_n$ is finite: 
$N=\#(\Gamma_n)+M$, where $\Gamma_n\subset\bZ^M$ is the set of
minimal solutions of the resonance equation {\rm(2.24)}.
\end{theorem}

Thus we conclude with the following.

The resonance algebra $\cF_n$ can be considered as a realization
of three different Poisson algebras:

--- the infinite Lie algebra (2.4),

--- the finite generated algebra (2.9) with a singular Poisson
tensor, 

--- the finite generated algebra (2.30) with constraints (2.29).

The first two Poisson algebras are universal and independent of
the frequency system~$n$. The third Poisson algebra depends
on~$n$ and has a smooth (polynomial) Poisson tensor.

\section{Resonance algebras for $2$-frequency systems}%3

In this section, we examine the $2$-frequency case. 
Besides the usual oscillator, we also consider the inverted
oscillator and the oscillator in a magnetic field.

\begin{example}[Resonance algebra of a two-dimensional oscillator]
In modern physics of nano-structures, 
the simplest model of a plane quantum dot 
or a plane artificial atom 
became very important because of the opening opportunities to
create these objects by new fine technologies.
The resonance case in this model is of special interest in view
of the effect of resonance stabilization (see \cite{n5}).
In this example 
we have the oscillator with two degrees of freedom:
\begin{equation}
H_0=\frac{n_1}{2}(q^2_1+p^2_1)+\frac{n_2}{2}(q^2_2+p^2_2),
\qquad
\wh{H}_0=n_1\hat{z}^*_1\hat{z}_1 +n_2\hat{z}^*_2\hat{z}_2,
\tag{3.1}
\end{equation}
where $n_1,n_2$ are coprime natural numbers. The resonance set
$\Gamma_n$ in this case consists of only two elements $\alpha$
and $-\alpha$, where $\alpha=(n_2,-n_1)\in\bZ^2$.
Thus the Poisson algebra $\cF_n$ is realized as the algebra of
functions on the constraint surface in $\bR^4$ given by 
Eq.~(2.29a): 
\begin{equation}
|\cA_\alpha|^2=\cA^{n_2}_1 \cA^{n_1}_2.
\tag{3.2}
\end{equation}
Here we use the simplified notation: 
$\cA_1\equiv \cA_{I_1}$, $\cA_2\equiv \cA_{I_2}$.
The Poisson brackets are 
\begin{align}
\{\cA_\alpha,\ol{\cA}_\alpha\}
&=i(n^2_2 \cA_2-n^2_1 \cA_1) \cA_1^{n_2-1} \cA_2^{n_1-1},
\nn\\
\{\cA_\alpha,\cA_1\}&=i n_2 \cA_\alpha,\qquad
\{\cA_\alpha,\cA_2\}=-i n_1 \cA_\alpha,
\tag{3.3}\\
\{\cA_1,\cA_2\}&=0.
\nn
\end{align}

In the general situation, for brackets like (2.30)
we can claim that 
the Jacobi condition holds on the constraint surface (2.29). 
But also, in the general situation, one can claim that 
it is sufficient to check the Jacobi condition only for a triple
of different indices from~$\Gamma_n$.
In our example, $\Gamma_n$ consists of two elements.
Thus the {\it Jacobi condition for the bracket {\rm(3.3)}
holds everywhere on~$\bR^4$}.
The constraint (3.2) in this case is actually the Casimir level.

There are two Casimir functions for the bracket (3.3):
\begin{equation}
C_0=n_1 \cA_1+n_2 \cA_2,\qquad
C_1=|\cA_\alpha|^2-\cA^{n_2}_1 \cA^{n_1}_2.
\tag{3.4}
\end{equation}
Note that the first function $C_0=\sum n_l \cA_l$ is a Casimir
function in the general case (2.30) as well.

The realization of the algebra (3.3) by functions
on $\bR^4$ is 
\begin{equation}
g_\alpha=z^{n_2}_1\oz^{n_1}_2,\qquad
g_{-\alpha}=\oz^{n_2}_1 z^{n_1}_2,\qquad
g_1=|z_1|^2,\qquad
g_2=|z_2|^2.
\tag{3.5}
\end{equation}
In this realization $C_1=0$ and $C_0=H_0$ is the Hamiltonian 
of the initial oscillator. 

The quantum version of generators (3.5) is 
\begin{equation}
\hat{g}_\alpha=\hat{z}^{*n_1}_2 \hat{z}^{n_1}_1,\qquad
\hat{g}_{-\alpha}=\hat{g}^{*}_\alpha,\qquad
\hat{g}_1=\hat{z}^{*}_1 \hat{z}_1,\qquad
\hat{g}_2=\hat{z}^{*}_2 \hat{z}_2.
\tag{3.6}
\end{equation}
They commute with the quantum oscillator (3.1).

The abstract version of these generators we denote by 
$\hat{\cA}_\alpha$, $\hat{\cA}^*_\alpha$, $\hat{\cA}_1$,
$\hat{\cA}_2$. 
They generate the quantum resonance algebra~$F_n$
of the two-dimensional oscillator. 

The commutation relations in this algebra are 
\begin{equation}
\begin{aligned}
{[}\hat{\cA}_\alpha,\hat{\cA}^*_\alpha]&=f(\hat{\cA}_1,\hat{\cA}_2),
\\
[\hat{\cA}_\alpha,\hat{\cA}_1]=\hb n_2\hat{\cA}_\alpha,
\qquad
[\hat{\cA}_\alpha,\hat{\cA}_2]&=-\hb n_1\hat{\cA}_\alpha,
\qquad
[\hat{\cA}_1,\hat{\cA}_2]=0.
\end{aligned}
\tag{3.7}
\end{equation}
Here the polynomial $f$ is defined by 
\begin{equation}
f(\cA_1,\cA_2)=\rho(\cA_1+\hb n_2,\,\cA_2-\hb n_1)
-\rho(\cA_1,\cA_2),
\tag{3.8}
\end{equation}
where
$$
\rho(\cA_1,\cA_2)\od \cA_1(\cA_1-\hb)\dots 
(\cA_1-(n_2-1)\hb)(\cA_2+\hb)\dots 
(\cA_2+n_1\hb).
$$
The Casimir elements of the algebra (3.7) are
\begin{equation}
\wh{C}_0=n_1\hat{\cA}_1+n_2\hat{\cA}_2,
\qquad
\wh{C}_1=\hat{\cA}^*_\alpha\hat{\cA}_\alpha 
-\rho(\hat{\cA}_1,\hat{\cA}_2).
\tag{3.9}
\end{equation}
In representation (3.6), 
$\hat{\cA}_\alpha\to\hat{g}_\alpha$, 
$\hat{\cA}_l\to\hat{g}_l$ the Casimir element $\wh{C}_1$ 
takes zero value.

Note that the polynomial $f$ on the right-hand side of relations
(3.7) has degree $n_1+n_2-1$.
In the case $n_1=n_2=1$ (i.e., 1:1 resonance) we obtain the Lie
algebra $\su(2)$.
In any other case the polynomial $f$ is of degree larger
than~$1$, and we have an algebra with nonlinear commutation
relations. 

The theory of irreducible representations and coherent states
for algebras of the type (3.7) was developed in \cite{kar7}.
Representations of the specific resonance algebra (3.7) 
in the case $n_1=1$, $n_2=2$ were described in \cite{n5}.
\end{example}

\begin{example}[Resonance algebra of the inverted oscillator]
As the next example, 
let us consider the inverted (unstable) oscillator
\begin{equation}
H^{(-)}_0=\frac12 \sum^{M}_{l=1} n_l (q^2_l - p^2_l).
\tag{3.10}
\end{equation}
This case is transformed to the previous one by the replacing 
$p\to ip$.
Thus we claim the set of functions in involution with
$H^{(-)}_0$ is given by the generators
\begin{align}
\cA^{(-)}_\alpha&\sim g^{(-)}_\alpha
\equiv\Big(\frac{q+p}{\sqrt{2}}\Big)^{\alpha_+}
\Big(\frac{q-p}{\sqrt{2}}\Big)^{\alpha_-},
\nn\\
\cA^{(-)}_{-\alpha}&\sim g^{(-)}_{-\alpha}
\equiv\Big(\frac{q+p}{\sqrt{2}}\Big)^{\alpha_-}
\Big(\frac{q-p}{\sqrt{2}}\Big)^{\alpha_+},
\tag{3.11}\\
\cA^{(-)}_l&\sim g^{(-)}_l=\frac{q^2_l-p^2_l}2,
\nn
\end{align}
where $n\circ\alpha=0$, $\alpha\in\bZ^M$. 
The commutation relations in the algebra $\cF_n$ are given by
a formula like (2.30) but without the imaginary unit~$i$ factor.
For instance, in the case $M=2$ we get the relation (3.3) without
the $i$-factor.

In order to understand the difference for which the $i$-factor
is responsible, let us consider the simplest 1:1 resonance.
Denote
\begin{equation}
\frac12(\cA_1-\cA_2)=b_3,\qquad
\frac1{2i}(\cA_\alpha-\cA_{-\alpha})=b_2,\qquad
\frac12(\cA_\alpha+\cA_{-\alpha})=b_1.
\tag{3.12}
\end{equation}
Then relations (3.3) read 
\begin{equation}
\{b_1,b_2\}=b_3,\qquad
\{b_2,b_3\}=b_1,\qquad
\{b_3,b_1\}=b_2.
\tag{3.13}
\end{equation}
The Casimir elements (3.4) are reduced to 
$C_1=b^2_1+b^2_2+b^2_3-C^2_0/4$.
In the realization (3.5), $C_1=0$, and so, 
$b^2_1+b^2_2+b^2_3=C^2_0/4$.

These are the commutation relations and Casimir elements 
of the Lie algebra $\su(2)$. 
Thus, in the case of the 1:1 resonance (for the usual
oscillator), the resonance algebra  is the envelope of the
$\su(2)$ Lie algebra. 

For the inverted oscillator in the 1:1 resonance case, 
we introduce the generators
\begin{equation}
b^{(-)}_1=\frac12(\cA^{(-)}_\alpha+\cA^{(-)}_\alpha),
\quad
b^{(-)}_2=\frac12(\cA^{(-)}_\alpha-\cA^{(-)}_{-\alpha}),
\quad
b^{(-)}_3=\frac12(\cA^{(-)}_1 - \cA^{(-)}_2),
\tag{3.14}
\end{equation}
where $\cA^{(-)}_\alpha$, $\cA^{(-)}_{-\alpha}$,
$\cA^{(-)}_1$, $\cA^{(-)}_2$ are given by (3.11). 
Then the relations between them are 
\begin{equation}
\{b^{(-)}_1,b^{(-)}_2\}=b^{(-)}_3,\qquad
\{b^{(-)}_2,b^{(-)}_3\}=b^{(-)}_1,\qquad
\{b^{(-)}_3,b^{(-)}_1\}=-b^{(-)}_2.
\tag{3.15}
\end{equation}
The Casimir elements in this case are reduced to 
$C_1=(b^{(-)}_1)^2 - (b^{(-)}_2)^2 + (b^{(-)}_3)^2 - C^2_0/4$.
In the realization (3.11), $C_1=0$, and so, 
$(b^{(-)}_1)^2 - (b^{(-)}_2)^2 + (b^{(-)}_3)^2 = C^2_0/4$.

The algebra (3.15) is the $\su(1,1)$ Lie algebra, and its
realization (3.11), (3.14) corresponds to the symplectic leaves
$\{(b^{(-)}_1)^2 - (b^{(-)}_2)^2 + (b^{(-)}_3)^2 = C^2_0/4\}$, 
which are one-sheeted hyperboloids.
\end{example}

\begin{example}
[Resonance algebra of 2D-oscillator in magnetic field]
To conclude this section, let us consider the $2$-dimensional
oscillator (the artificial atom) placed 
into a magnetic field. 
Let us assume that the magnetic field is homogeneous and
perpendicular to the plane where the oscillator sits.
Denote by $\omega_L$ the half of the Larmor frequency (which is
the magnetic field strength multiplied by the charge of the
particle and divided by its mass). 
Let the oscillator be isotropic, and let $\omega_0$ be its 
frequency. If we represent the ratio $\omega_L/\omega_0$ as 
\begin{equation}
\Big(\frac{\omega_L}{\omega_0}\Big)^2
=\frac{s^2}{k^2-s^2},\qquad k>s>0,
\tag{3.16}
\end{equation}
then the Hamiltonian of this physical model reads
\begin{equation}
H_0=k\Big(\frac{|q|^2+|p|^2}{2}\Big)+s(q_1p_2-q_2p_1)
=k(|z_1|^2+|z_2|^2)+is(\oz_2 z_1 - \oz_1 z_2).
\tag{3.17}
\end{equation}

A resonance happens if the numbers $k,s$ are commensurable, and
so one can find a pair of coprime numbers $l,m$ such that 
\begin{equation}
\frac{k+s}{k-s}=\frac{l}{m}.
\tag{3.18}
\end{equation}

Note that in the coordinates 
$z_{\pm}=\frac1{\sqrt{2}}(z_1\mp i z_2)$ 
the system is equivalent to the oscillator 
$l|z_+|^2+m|z_-|^2$.
Its resonance algebra was already described above
(in Example~3.1).
Thus we can claim that the {\it resonance algebra $\cF_{l,m}$ of
the Hamiltonian\/} (3.17) 
{\it under condition\/} (3.18) {\it has the generators\/}:  
\begin{equation}
\cA_1=|z_+|^2,\qquad 
\cA_2=|z_-|^2,\qquad 
\cA_+=\oz^l_- z^m_+,\qquad 
\cA_-=\ol{\cA}_+.
\tag{3.19}
\end{equation}
The brackets between them are 
\begin{equation}
\begin{aligned}
\{\cA_+,\cA_1\}&=im \cA_+,\qquad 
\{\cA_+,\cA_2\}=-il \cA_+,\qquad 
\{\cA_1,\cA_2\}=0,
\\
\{\cA_+,\cA_-\}&=i(m^2 \cA_2 - l^2 \cA_1) \cA^{m-1}_1 \cA^{l-1}_2.
\end{aligned}
\tag{3.20}
\end{equation}

This algebra just coincides with the algebra (3.3) if $n_1=l$,
$n_2=m$. Thus it has the same Casimir elements as in (3.4):
$$
C_0=l \cA_1 + m \cA_2,\qquad
C_1 = \cA_+ \cA_- - \cA^m_1 \cA^l_2,
$$
and $C_1=0$ on the realization (3.19).

The quantum version of the algebra (3.20) is easily obtained 
in the same way as in Example~3.1 (see (3.7)).

Note that the nonlinear character of the Poisson tensor in
(3.20) appears automatically if $l\neq m$ or $s\neq 0$ in
(3.17), i.e., if the magnetic field is not zero. 
The magnetic field is responsible for this nonlinearity, 
in spite of the rotational symmetry of the system (3.17).

We can also consider a non-isotropic oscillator 
with two different frequencies placed into a magnetic field:
\begin{equation}
H_0=\frac12|p|^2 + \frac12(\omega^2_1 q^2_1 + \omega^2_2 q^2_2) 
+ (q_1 p_2 - q_2 p_1).
\tag{3.21}
\end{equation}
The effective frequencies $\omega_{\pm}$ of the Hamilton flow in
this model are given by the formula
\begin{equation}
\omega^2_{\pm}=\frac12\Big[\omega^2_1 + \omega^2_2 +2 
\pm\sqrt{(\omega^2_1 - \omega^2_2)^2 
+ 8(\omega^2_1 + \omega^2_2)}\,\Big].
\tag{3.22}
\end{equation}
The resonance happens if $\omega_+/\omega_-=l/m$, 
where $l,m$ are certain coprime numbers. One can describe the
resonance algebra of this nonsymmetric system in the same way as
above.
\end{example} 

\section{Resonance precession}%4

Now let us discuss some spectral and dynamical effects
accompanying the appearance of a resonance. 

For simplicity, we consider only the case of two degrees of
freedom. Let a Hamiltonian $H$ have the nondegenerate minimum at
the point $x=0$. Around this point, we have the Taylor expansion:
\begin{equation}
H(x)=\const +H_0(x)+H_1(x)+\dots,
\tag{4.1}
\end{equation}
where $H_j$ are the homogeneous polynomials of degree $j+2$.
By passing to the new variables
$$
x=\hb^{1/N} x',\qquad \hb'=\hb^{1-2/N}\qquad (N\geq2),
$$
we obtain 
\begin{equation}
H=\const + \hb^{2/N}\big(H_0(x')
+\hb^{1/N} H_1(x') +\hb^{2/N} H_2(x') +\dots\big).
\tag{4.2}
\end{equation}

Let the variation matrix of the Hamilton system at the
equilibrium point $x'=0$ have only imaginary eigenvalues 
$\pm i\omega_1$, $\pm i\omega_2$. 
Then in suitable coordinates one can represent the quadratic
form $H_0$ as the oscillator
\begin{equation}
H_0=\omega_1|z'_1|^2 + \omega_2|z'_2|^2,
\tag{4.3}
\end{equation}
where $z'_1$, $z'_2$ are complex coordinates of the point 
$x'\in\bR^4$.

If the frequencies $\omega_1,\,\omega_2$ are not in resonance,
then, applying the averaging method in the $N$th microzone, 
where $x=O(\hb^{1/N})$, 
one can reduce the Hamiltonian to a function in the action
variables $\cA_1=|z'_1|^2$ and $\cA_2=|z'_2|^2$ only:
\begin{equation}
H\sim\const+\hb^{2/N}H_0+\hb^{3/N} f_1(\cA_1,\cA_2)
+\hb^{4/N} f_2(\cA_1,\cA_2)+\dots\,.
\tag{4.4}
\end{equation}
The same is true on the quantum level.

The action variables $\cA_1,\cA_2$ in (4.4) commute with each other
and thus the dynamics and the algebra are reduced to purely
commutative ones.

Let us assume now that the frequencies $\omega_1$, $\omega_2$ are
in a resonance 
\begin{equation}
\frac{\omega_1}{\omega_2}=\frac{n_1}{n_2},\qquad
\text{$n_1$ and $n_2$ are coprime integers}.
\tag{4.5}
\end{equation}
Then the application of the averaging method to (4.2) implies
(see in Appendix) that the Hamiltonian is reduced to 
\begin{equation}
H\sim\const+\hb^{2/N} H_0
+\hb^{3/N} f(\cA_1,\cA_2,\cA_\alpha,\ol{\cA}_\alpha),
\tag{4.6}
\end{equation}
where $f=f_1+\hb^{1/N}f_2+\dots$ and 
$f_j$ are functions in generators of the
resonance algebra $\cF_n$ with relations (3.3). 
The same representation exists in the quantum case as well, and
in this case, one has to use the resonance algebra (3.7)
(with the parameter $\hb'$ instead of~$\hb$).

Thus under the resonance condition (4.5), 
instead of the commutative algebra of action variables, 
we have the noncommutative algebra (3.3).
The averaged system (see the definition in \cite{kar15}) 
can be represented as a Hamiltonian system 
on symplectic leaves
\begin{equation}
\Omega=\{C_0=\const, C_1=0\},
\tag{4.7}
\end{equation}
where $C_0$, $C_1$ are the Casimir funcitons (3.4).
These leaves are compact surfaces in $\bR^4$ diffeomorphic to
the sphere $\bS^2$.  

The Hamiltonian system corresponding to the perturbing
term $f$ in (4.6) reads
\begin{equation}
\frac{d}{dt}\cA=\{f,\cA\}
\tag{4.8}
\end{equation}
or explicitly,
\begin{align}
\frac{d}{dt}\cA_\alpha
&=in_1\frac{\pa f}{\pa \cA_2}-in_2\frac{\pa f}{\pa \cA_1}
+i(n^2_1 \cA_1 - n^2_2 \cA_2) \cA^{n_2-1}_1 \cA^{n_1-1}_2
\frac{\pa f}{\pa\ol{\cA}_\alpha},
\nn\\
\frac{d}{dt}\cA_1
&=in_2\cA_\alpha\frac{\pa f}{\pa \cA_\alpha}
-in_2\ol{\cA}_\alpha\frac{\pa f}{\pa\ol{\cA}_\alpha},
\tag {4.8a}\\
\frac{d}{dt}\cA_2
&=-in_1\cA_\alpha\frac{\pa f}{\pa \cA_\alpha}
+in_1\ol{\cA}_\alpha\frac{\pa f}{\pa\ol{\cA}_\alpha}.
\nn
\end{align}

In the particular case $n_1=n_2=1$ (the isotropic resonance 1:1)
the algebra (3.3) is reduced to (3.13), i.e., to $\su(2)$. 
Then the leaves (4.7) are exactly the spheres and the system
(4.8) is the classical Euler system for the spinning top 
\cite{n7}.

In the case of general anisotropic resonance $n_1:n_2$, 
we obtain something like a generalized top with non-Lie
commutation relations (3.3) and the Hamilton dynamics (4.8) 
on it. 
{\it This dynamics describes the evolution of the noncommutative
integrals of motion of the resonance oscillator}. 

In the phase space $\bR^4_x$ on each energy level 
$H_0\equiv C_0=\const$, we have a fibration by periodic
trajectories of $H_0$. 
The whole variety of trajectories is parametrized by the
coordinates $\cA_1,\cA_2,\cA_\alpha$. 
The evolution (4.8) is similar to the evolution 
of the rotation axis in the theory of spinning top.
That rotation is known as precession.
Thus in our case, where we have a resonance analog of the
spinning top, one can use the term {\it resonance precession\/}
for solutions of the evolution system (4.8).

The quantum version of (4.8) is the Heisenberg type equations:
\begin{equation}
-i\hb'\frac{d}{dt}\hat{\cA}(t) = [\hat{f},\hat{\cA}(t)],
\qquad
\hat{\cA}(t)\big|_{t=0}=\hat{\cA}.
\tag{4.9}
\end{equation}
Here $\hat{\cA}(t)$ denotes the evolution of generators
$\hat{\cA}$ of the algebra (3.7). The Hamiltonian $\hat{f}$ in
(4.9) is a polynomial function in generators of the same algebra. 
The system (4.9) describes the {\it quantum resonance precession}.

Equations (4.9) are related to the $N$th microzone with $N>2$.
If $N=2$ (this is called the nanozone, \cite{kar21,n5}), 
then $\hb'=1$ and (4.9) has to be read as
\begin{equation}
-i\frac{d}{dt}\hat{\cA}(t)
=[\hat{f},\hat{\cA}(t)].
\tag{4.9a}
\end{equation}

Studying the spectral problem for the Hamiltonian $\wh{H}$, we
can exploit the representation (4.6) as was demonstrated in
\cite{kar21,n5}. In this situation, instead of the dynamical system
(4.9), (4.9a),  
one has to consider the spectral problem for the operator 
$\hat{f}$ over the resonance algebra. This problem can
be effectively analyzed by applying the coherent states method,
see \cite{kar21,n5}. 

Note that there are some simple geometrical and topological
implications of the compactness of symplectic leaves $\Omega$,
where the resonance precession system (4.8) is actually
supported. For instance, on $\Omega$, the Hamiltonian $f$ must
have stationary points (at least two of them), which provide the
equilibria in the $a$-space. 
The quantum operator $\hat{f}$ near these points can be
analyzed using the method of deformed coherent states
\cite{kar20}, which generalizes the well-known complex germ
method \cite{kar16} to the case of algebras with non-Heisenberg
relations like~(3.7). 

\section[Triple algebras and solving the resonance precession\\
system via noncommuting space--time coordinates]{Triple algebras
and solving the resonance precession 
system via noncommuting\\ space--time coordinates}%5

Now we investigate the resonance precession in more detail.
First of all, let us analyze the averaging transformation which
produces the precession Hamiltonian over the resonance algebra
from a Hamiltonian over the original (Heisenberg) algebra. 

Let the leading part of the expansion (4.1) be the resonance
oscillator  
$H_0=\sum^{M}_{l=1} n_l|z_l|^2$.
Then the averaging transformation (denoted by $\Pi$) is given by
formula (A.7) of the Appendix:
\begin{equation}
F=\sum_k c_k g_k\quad \overset{\Pi}{\to}\quad 
\uF =\sum_{k\in R_n} c_kg_k.
\tag{5.1}
\end{equation}
Here $c_k\in \bC$, the monomials $g_k(z)=z^{k_+}\oz^{k_-}$ 
are determined by multi-indices 
$k=(k_+,k_-)\in\bZ^M_+\times\bZ^M_+$,
and the resonance set $R_n$ is determined by solutions of the
resonance equation (1.10) $n\circ(k_+ - k_-)=0$.

Denote by $\cF$ the Poisson algebra of all (complex) polynomials
over $\bR^{2M}$. It is split into the direct sum of subspaces 
$$
\cF=\cL\oplus\cL^\perp,\qquad
\cL=\Pi(\cF),\qquad
\cL^\perp=(I-\Pi)(\cF).
$$
Thus the averaging transformation (5.1) is the projection of
$\cF$ onto $\cL$. 

The multiplication operation in $\cF$ and the Poisson brackets
are naturally projected into $\cL$: 
\begin{equation}
\begin{aligned}
F&=\sum_{k\in R_n}c_k g_k, \qquad
G=\sum_{k\in R_n} d_k g_k \quad
\Longrightarrow 
\\
FG&=\sum_{k,r\in R_n} c_k d_r g_{k+r},\qquad
\{F,G\}_\cL \od\sum_{k,r\in R_n} c_k d_r\{g_k,g_r\}_\cF.
\end{aligned}
\tag{5.2}
\end{equation}

The general mechanism of such a projection of the Poisson
algebra structure can be easily described.

\begin{lemma}
Let $\cF$ be an abstract Poisson algebra which admits the
splitting into the direct sum of linear subspaces
$$
\cF=\cL\oplus\cL^\perp, 
$$
and let $\cL^\perp$ be the Poisson $\cL$-modulus
\begin{equation}
\cL^\perp\circ\cL\subset\cL^\perp,\qquad
\{\cL^\perp,\cL\}\subset\cL^\perp.
\tag{5.3}
\end{equation}
Then there is a unique Poisson algebra structure on $\cL$ such
that  
\begin{equation}
F\underset{\cL}{\circ}G=\Pi(F\underset{\cF}{\circ}G),
\qquad
\{F,G\}_\cL=\Pi(\{F,G\}_\cF),\qquad\forall F,G\in\cL,
\tag{5.4}
\end{equation}
where $\Pi:\cF\to\cL$ is the projection on $\cL$ along 
$\cL^\perp$.

One has the identity 
\begin{equation}
\begin{aligned}
\Pi(F\underset{\cF}{\circ}G)
&=\Pi(F)\underset{\cL}{\circ}\Pi(G)
+\Pi\big( (I-\Pi)(F)\underset{\cF}{\circ}(I-\Pi)(G)\big),
\\
\Pi(\{F,G\}_\cF)
&=\{\Pi(F),\Pi(G)\}_\cL
+\Pi\big(\{(I-\Pi)(F),(I-\Pi)(G)\}_\cF\big).
\end{aligned}
\tag{5.5}
\end{equation}
If $\cL$ is a Poisson subalgebra in $\cF$, then 
$$
F\underset{\cL}{\circ}G=F\underset{\cF}{\circ}G,\qquad
\{F,G\}_\cL=\{F,G\}_\cF,\qquad\forall F,G\in\cL.
$$
\end{lemma}

Indeed, one can take (5.4) just as the definition of the
multiplication operation and the bracket on $\cL$, 
and check that condition (5.3) implies the associativity and the
Jacobi identity for this multiplication and the bracket.

In the case of the averaging projection (5.1), 
the multiplication and the bracket (5.4)
are given by formula (5.2).
The identities (5.5) read
\begin{equation}
\begin{aligned}
\Pi(FG)&=\Pi(F)\Pi(G)
+\sum_{k,r\not\in R_n, k+r\in R_n} c_k d_r g_{k+r},
\\
\Pi(\{F,G\}_\cF)&=\{\Pi(F),\Pi(G)\}_\cL
+\sum_{k,r\not\in R_n, k+r\in R_n} c_k d_r 
\sum^{M}_{l=1} [k,r]_l\, g_{k+r-I_l}.
\end{aligned}
\tag{5.5a}
\end{equation}
The last summands in the right-hand sides of (5.5a) are the 
{\it anomaly\/} which prevents the projection $\Pi:\cF\to\cL$
from being a homomorphism. 
In particular, this anomaly makes the Poisson bracket in $\cL$  
be noncommutative in the resonance case (but the associative
multiplication in $\cL$ is, of course, Abelian).

In general, the projection allows one to define nontrivial
brackets on linear subspaces in $\cF$.

\begin{lemma}
Let $\cL$ obey condition {\rm(5.3)} and
$\cF_0\subset\cF$ be a linear subspace
such that the projection $\Pi:\cF\to\cL$ is injective on $\cF_0$
and the image $\cL_0\od \Pi(\cF_0)$ is a Poisson subalgebra in
$\cL$.  Then $\cF_0$ is endowed with a Poisson algebra structure
which makes $\Pi:\cF_0\to\cL_0$ be a homomorphism. 
\end{lemma}

This trivial statement can be amplified as follows.

\begin{lemma}
Let a subspace $\cL\subset\cF$ obey {\rm(5.3)}, 
$\Pi:\cF\to\cL$ be injective on a subspace
$\cF_0\subset\cF$, 
and $\cL$ be split into the direct sum of
linear subspaces,
\begin{equation}
\cL=\cL_0\oplus\cL_1
\tag{5.6}
\end{equation}
such that $\cL_0=\Pi(\cF_0)$ is the subalgebra with respect to
the multiplication in $\cL$, and 
\begin{equation}
\{\cL_0,\cL_1\}_\cL\subset\cL_0.
\tag{5.7}
\end{equation}
Then on $\cL_0$, and so on the subspace $\cF_0$, 
there exists the triple Poisson algebra structure with the
following triple bracket 
\begin{equation}
\{F,G,E\}_{\cL_0}\od \{ \{F,G\}_\cL,E\}_\cL,\qquad 
F,G,E\in\cL_0.
\tag{5.8}
\end{equation}
\end{lemma}

In this framework, let $\cF$ be the algebra of polynomial
functions on the phase space $\bR^{2M}=\bR^M_q\oplus\bR^M_p$, 
the projection $\Pi$ be defined by the averaging transformation
(5.1), and $\cF_0$ be a subspace consisting of functions
constant along $\bR^M_p$. 
Then $\cL=\cF_n$ is the resonance algebra. If there is a
splitting of the resonance algebra like (5.6), (5.7), then the
triple bracket (5.8) generates a noncommutative
triple algebra structure on $\cF_0$. 
This might be interpreted as the appearance of a triple
{\it noncommutative structure on the configuration space\/} 
$\bR^M_q$ (since $\cF_0$ consists of functions in
$q$-coordinates). 

In the simple situation, where we know that the averaging of
$\cF_0$-functions forms a Poisson subalgebra in the resonance
algebra $\cF_n$, one can apply Lemma~5.2. 
Then we can claim that the space $\cF_0$ of functions on
$\bR^M_q$ is a Poisson noncommutative algebra, and so, 
the configuration space is noncommutative in the usual sense. 
The well-known example of coordinate noncommutativity of this
type is presented by the Landau model (a charged particle on a
plane with a perpendicular magnetic field \cite{n11,GJ,J}).
In this example, 
we have the case of a neutral resonance, where the variation
matrix of the Hamiltonian system related to $H_0$ 
has the twice degenerate zero eigenvalue, see \cite{n5}.

The resonance considered in the given paper belongs to
the stable (elliptic) class: the eigenvalues $\pm in_l$ of the 
variation matrix are imaginary nonzero numbers. 
In this situation, Lemma~5.2 cannot be applied any more, but
Lemma~5.3 works.
We describe now simple examples of the 1:1 and 1:2 resonances,
demonstrate the triple brackets over the configuration plane
$\bR^2_q$, 
and show how the resonance precession dynamics can be resolved
on this noncommutative plane.

\begin{example}[Noncommutative space--time under 1:1 resonance] 
Let us consider the isotropic resonance oscillator 
$H_0=\frac12(q^2_1+p^2_1)+\frac12(q^2_2+p^2_2)$.
The averaging transformation $\Pi$ acts via the integral
\begin{equation}
\Pi(F)(q,p)=\frac1{2\pi}\int^{2\pi}_{0}
F(q\cos t+p\sin t,p\cos t-q\sin t)\,dt.
\tag{5.9}
\end{equation}
So, the average of any odd function $F(x)=-F(x)$ is zero:
$\Pi(F)=0$.
Therefore, we consider the subspace 
$\cF_0\subset\cF(\bR^2_q\times \bR^2_p)$ consisting of all even
functions on $\bR^2_q$. This subspace is generated by 
$q^2_1$, $q^2_2$, $q_1q_2$.
Their averages we denote by
\begin{equation}
\begin{aligned}
X&=\underline{q^2_1}=\frac12(q^2_1+p^2_1),\qquad
Y=\underline{q^2_2}=\frac12(q^2_2+p^2_2),\\
Z&=\underline{q_1q_2}=\frac12(q_1q_2+p_1p_2).
\end{aligned}
\tag{5.10}
\end{equation}
Also introduce\footnote{The previous
notation $X\equiv \cA_1$, $Y\equiv \cA_2$, 
$Z\equiv\frac12(\cA_\alpha+\ol{\cA}_\alpha)$, and 
$W\equiv\frac1{2i}(\cA_\alpha-\ol{\cA}_\alpha)$ 
see in Section~3.}  
the function $W=\frac12(p_1q_2-q_1p_2)$. 
All together $X,\,Y,\,Z,\,W$ generate the
subalgebra $\cL\subset \cF$, this is exactly the resonance
algebra $\cL\equiv\cF_{1,1}$ of the isotropic oscillator $H_0$
with frequencies 1:1. 
Now, $\cL$ can be split into the sum of the subspace $\cL_0$
enveloping $X,Y,Z$ and the subspace $\cL_1$ enveloping $W$.

We can check condition (5.7) and evaluate the triple bracket
(5.8). 
Indeed, the brackets between generators of the resonance algebra
are
\begin{equation}
\begin{aligned}
\{X,Y\}&=0,&\qquad\{X,Z\}&=W,&\qquad\{Y,Z\}&=-W,
\\
\{X,W\}&=-Z,&\qquad\{Y,W\}&=Z,&\qquad\{Z,W\}&=\frac12(X-Y).
\end{aligned}
\tag{5.11}
\end{equation}
The Casimir functions are 
$$
C_0=X+Y,\qquad C_1=XY-Z^2-W^2.
$$
From (5.11) one obtains the triple brackets
\begin{equation}
\begin{aligned}
\{\{X,Z\},X\}&=Z,&\qquad\{\{X,Z\},Y\}&=-Z,&\qquad\{\{X,Z\},Z\}&=\frac12(Y-X),
\\
\{\{Y,Z\},X\}&=-Z,&\qquad\{\{Y,Z\},Y\}&=Z,&\qquad\{\{Y,Z\},Z\}&=\frac12(X-Y).
\end{aligned}
\tag{5.12}
\end{equation}

By passing from (5.11) to (5.12) we exclude the momentum
coordinate~$W$.

Thus, over the configuration plane $\bR^2_q$ 
we have the triple Poisson algebra $\cF_0$ (5.12) with the
Casimir function $C_0=X+Y$. 
The noncommuting coordinates $X,Y,Z$ on $\bR^2_q$ 
are the mean square deviations of the Euclidean coordinates.

Let us write the resonance precession system in these
coordinates. 
Assume that the perturbing Hamiltonians in (4.2) are all coming
from the expansion of a potential well~$V$ and so they are 
$p$-independent:
$$
H_j(x)=\sum_{|\alpha|=j+2}\frac1{\alpha!} D^\alpha V(0) q^\alpha.
$$
This means that the total Hamiltonian $H$ is
\begin{equation}
H=H_0+V=\frac12(|q|^2+|p|^2)+V(q),\qquad q,p\in\bR^2,
\tag{5.13}
\end{equation}
where the potential $V$ contains terms of the third and higher
orders near the point $q=0$.

The averaging transformation (5.9) being applied to $H_1, H_3,\dots$
gives zero. The leading nonzero term is
$$
\underline{H_2}\equiv \Pi(H_2)
=\sum_{|\alpha|=4}\frac1{\alpha!} D^\alpha V(0) \underline{q^\alpha}.
$$ 
The averages $\underline{q^\alpha}$ are given in the following list:
\begin{equation}
\begin{gathered}
\underline{q^4_1}=\frac32 X^2,\qquad 
\underline{q^4_2}=\frac32 Y^2,
\\
\underline{q_1 q^3_2}=\frac32 YZ,\qquad
\underline{q^3_1 q_2}=\frac32 XZ,\qquad
\underline{q^2_1 q^2_2}=\frac12 XY+Z^2.
\end{gathered}
\tag{5.14}
\end{equation}
It is important to mention that there is no generator $W$ in
these formulas. 

Finally, the average of the perturbing Hamiltonian reads
\begin{equation}
f(X,Y,Z)=\alpha X^2+\beta Y^2+\gamma Z^2+\frac12\gamma XY
+\delta XZ+\rho YZ+O^4,
\tag{5.15}
\end{equation}
where 
\begin{align*}
\alpha&=\frac1{16}\frac{\pa^4 V}{\pa q^4_1}(0),&\qquad
\beta&=\frac1{16}\frac{\pa^4 V}{\pa q^4_2}(0),&\qquad
\gamma&=\frac1{4}\frac{\pa^4 V}{\pa q^2_1\pa q^2_2}(0),
\\
\delta&=\frac1{4}\frac{\pa^4 V}{\pa q^3_1\pa q_2}(0),&\qquad
\rho&=\frac1{4}\frac{\pa^4 V}{\pa q_1\pa q^3_2}(0),
\end{align*}
and $O^4$ denotes terms of order~$4$ and higher near the
origin~$0$. 

Now we can write the resonance precession system (4.8) for the
coordinates $X,Y,Z$, and $W$. 
Actually, 
because of the Casimir constraints $C_1=0$ and $X+Y=C_0$,
it is enough to consider the equation for coordinates 
$a\od X-Y$ and $b\od 2Z$ only. 
These coordinates characterize the shape of the oscillator 
orbits projected to $\bR^2_q$ (the eccentricity and the shear).

The resonance precession system reads
\begin{equation}
\frac{d}{dt}a=-4W\frac{\pa f}{\pa b},\qquad
\frac{d}{dt}b=4W\frac{\pa f}{\pa a},\qquad
\frac{dW}{dt}=a\frac{\pa f}{\pa b}-b\frac{\pa f}{\pa a}.
\tag{5.16}
\end{equation}
Here $f$ is taken from (5.15) and expressed 
in the $a,b$-coordinates.

In order to exclude $W$, 
let us introduce a new ``time'' $\tau$ on the space--time 
$\bR^2_{a,b}\times\bR_t$ by means of the equation
$$
\{f,\tau\}=-4W.
$$
Taking the second bracket, we obtain 
\begin{equation}
\{\{f,\tau\},f\}=4\bigg(a\frac{\pa f}{\pa b}-b\frac{\pa f}{\pa a}\bigg).
\tag{5.17}
\end{equation}
Note that the new time $\tau$ does not commute with $a,b$ and we have
\begin{equation}
\begin{aligned}
\{\{f,\tau\},a\}&=-4b,&\qquad\{\{f,\tau\},b\}&=4a,
\\
\{\{a,b\},b\}&=-4a,&\qquad\{\{a,b\},a\}&=4b.
\end{aligned}
\tag{5.18}
\end{equation}

On the trajectories of (5.16), relation (5.17) becomes 
an equation for~$\tau$
\begin{equation}
\frac{d^2\tau}{dt^2}=4\bigg(b\frac{\pa f}{\pa a}
-a\frac{\pa f}{\pa b}\bigg)\bigg|_{\substack{a=a(\tau),\\ b=b(\tau)}}.
\tag{5.19}
\end{equation}
The trajectory $(a(\tau),b(\tau))$ parametrized by the new time
is obtained from the first two equations in (5.16):
\begin{equation}
\frac{da}{d\tau}=\frac{\pa f}{\pa b},
\qquad 
\frac{db}{d\tau}=-\frac{\pa f}{\pa a}.
\tag{5.20}
\end{equation}
Using (5.20) we resolve (5.19) as follows:
\begin{equation}
\frac12\int^{\tau}_{0}
\frac{d\tau}{\sqrt{C^2_0-a(\tau)^2-b(\tau)^2}}=t.
\tag{5.21}
\end{equation}

Now we note that the function $f$ (5.15) after the change of
variables $X=\frac12(C_0+a)$, $Y=\frac12(C_0-a)$, $Z=\frac12b$ 
is quadratic in $a,b$ (in the classical limit in any $N$th
microzone near the origin).
Thus, the right-hand side of (5.20) is just linear in $a,b$ 
and so system (5.20) is explicitly integrable in trigonometric
functions of time~$\tau$. This 
{\it completely resolves the resonance
precession system for the Hamiltonian}~(5.13).

Note that the reduced system (5.20) can be considered as the
Hamiltonian system of the configuration space--time 
with the triple Poisson structure (5.18).
\end{example}

\begin{example}[Noncommutative space--time under resonance 1:2]
Let us consider the oscillator 
$H_0=\frac12(q^2_1+p^2_1)+(q^2_2+p^2_2)$.
Applying the averaging transformation (5.1), we obtain the
following coordinates on the configuration plane:
$$
X=\underline{q^2_1},\qquad 
Y=\underline{q^2_2},\qquad 
Z=\underline{q^2_1 q_2}.
$$
The explicit formulas for them via the coordinates on 
$\bR^2_q\times\bR^2_p$ are
$$
X=\frac12(q^2_1+p^2_1),\qquad
Y=\frac12(q^2_2+p^2_2),\qquad
Z=\frac14(q^2_1q_2+2q_1p_1p_2-q_2p^2_1).
$$
We also introduce 
$W=\frac14(p^2_1p_2+2q_1q_2p_1-q^2_1p_2)$.
Then the mutual brackets between these generators 
of the resonance algebra $\cF_{1,2}$ are the following:
\begin{equation}
\begin{aligned}
\{X,Y\}&=0,\qquad\{X,Z\}=2W,\qquad\{Y,Z\}=-W,
\\
\{X,W\}&=-2Z,\qquad\{Y,W\}=Z,\qquad\{Z,W\}=\frac14X^2-XY.
\end{aligned}
\tag{5.22}
\end{equation}
The Casimir functions for this quadratic Poisson bracket are 
$C_0=X+2Y$ and $C_1=\frac12 X^2Y-Z^2-W^2$.

By excluding $W$ we obtain the triple brackets
\begin{align}
\{\{X,Z\},X\}&=4Z,&\qquad \{\{X,Z\},Y\}&=-2Z,
\nn\\
\{\{Y,Z\},X\}&=-2Z,&\qquad \{\{Y,Z\},Y\}&=Z,
\tag{5.23}\\
\{\{X,Z\},Z\}&=2XY-\frac12X^2,&\qquad \{\{Y,Z\},Z\}&=\frac14X^2-XY.
%\nn
\end{align}
Note that these triple brackets are nonlinear (quadratic) in
coordinates. So, this is not a Lie triple system, but its
nonlinear generalization.

The resonance precession system in this case also has the form
(5.20), where the coordinates $a,b$ are defined by 
$$
a=X-2Y,\qquad b=Z.
$$

The new time $\tau$ in (5.20) is determined by the equations
$$
\frac{d\tau}{dt}=-4W,\qquad 
\frac{dW}{dt}=\bigg(\frac14X^2-XY\bigg)\frac{\pa f}{\pa Z}
-Z\bigg(2\frac{\pa f}{\pa X}-\frac{\pa f}{\pa Y}\bigg),
$$
which are integrated similarly as in (5.21).

In conclusion, we remark that the noncommuting coordinates and
the triple brackets which we introduced on the configuration
plane $\bR^2_q$ admit natural quantum analogs. It is easy to
derive quantum versions of (5.12) and (5.23).
\end{example}

\addcontentsline{toc}{section}{Appendix. Operator averaging}

\section*{Appendix. Operator averaging}

Let us consider the family of operators 
$\wh{H}_0+\ve\wh{H}_1$, $\ve\to0$, 
in the algebra~$F$.
We are looking for a family of operators 
$\wh{U}_\ve$ such that 
\begin{equation}
(\wh{H}_0+\ve\wh{H}_1)\wh{U}_\ve
=\wh{U}_\ve(\wh{H}_0+\ve\underline{\wh{H}_1}
+\ve^2 \underline{\wh{H}_2}+O(\ve^3)),
\tag{A.1}
\end{equation}
where
\begin{equation}
[\wh{H}_0,\underline{\wh{H}_1}]=0,\qquad
[\wh{H}_0,\underline{\wh{H}_2}]=0.
\tag{A.2}
\end{equation}
(we stopped at the order $O(\ve^3)$ only for simplicity).
One can formally choose 
$\wh{U}_\ve=\exp\{-i\ve(\hat{f}_0+\ve\hat{f}_1+O(\ve^2))\}$. 
In more explicit form, 
\begin{equation}
\wh{U}_\ve=I-i\ve\hat{f}_0
-\ve^2(i\hat{f}_1+\hat{f}^2_0/2)+O(\ve^3), 
\tag{A.3}
\end{equation}
where $\hat{f}_0,\hat{f}_1$ are solutions of the ``homological''
equations: 
\begin{equation}
i[\wh{H}_0,\hat{f}_0]=\wh{H}_1-\underline{\wh{H}_1},
\qquad
i[\wh{H}_0,\hat{f}_1]=\wh{H}_2-\underline{\wh{H}_2}.
\tag{A.4}
\end{equation}
here we denote $\wh{H}_2\od\frac{i}2[\hat{f}_0,
\wh{H}_1+\underline{\wh{H}_1}]$.

Thus one has to solve Eqs.~(A.4) with additional conditions (A.2).

Assume now that $\wh{H}_0$ is a quadratic form: 
$\wh{H}_0=\sum^{M}_{l=1}\omega_l\hat{z}^*_l\hat{z}_l$
(see (1.1a), (1.2a)) with positive frequencies $\omega_l$.

Any operator from the algebra $F$, say $\wh{H}_1$,
can be represented as 
\begin{equation}
\wh{H}_1=\sum_k c_k \hat{g}_k,\qquad
\hat{g}_k=\hat{z}^{*k_+}\hat{z}^{k_-},
\tag{A.5}
\end{equation}
where $c_k\in\bC$, $k_+,k_-\in\bZ^M_+$. 
Then one solves Eqs.~(A.4), (A.2) by the explicit formulas
\begin{equation}
\hat{f}_0=\frac1{i}\sum_{\omega\circ(k_+-k_-)\neq0}
\frac{c_k}{\omega\circ(k_+-k_-)}\hat{g}_k,
\qquad
\underline{\wh{H}_1}=\sum_{\omega\circ(k_+-k_-)=0}
c_k \hat{g}_k.
\tag{A.6}
\end{equation}
The second equations in (A.4), (A.2) are solved similarly.

The transformation (A.1) is called the {\it operator averaging}.
If the operator $\wh{U}_\ve$ (A.3) is invertible, one can
conclude 
from (A.1) that the original operator $\wh{H}_0+\ve\wh{H}_1$ 
is equivalent (up to $O(\ve^3)$) to the new one 
$\wh{H}_0+\ve\underline{\wh{H}_1}+\ve\underline{\wh{H}_2}$, 
which contains the perturbations $\underline{\wh{H}_1}$
and $\underline{\wh{H}_2}$ commuting with the leading part 
$\wh{H}_0$. 
Thus, $\underline{\wh{H}_1}$ and $\underline{\wh{H}_2}$ are 
elements of the commutant $F_\omega$. 

From Sections~1 and~2 we know how this commutant looks like 
(see Theorem~1.7) and can represent 
$\underline{\wh{H}_1}$, $\underline{\wh{H}_2}$ as functions 
in generators of its resonance components 
$F^{(1)},\dots,F^{(L)}$.
This representation follows from (A.6) and (1.7)
\begin{equation}
\underline{\wh{H}_1}
=\sum_{\substack{k^{(j)}\in R_{n^{(j)}},\\ j=1,\dots,L}} 
c_k\hat{g}_k,
\tag{A.7}
\end{equation}
where $k=(k^{(1)},\dots,k^{(L)})$ 
and $R_{n^{(j)}}$ is the resonance set (1.10) corresponding to
the $j$th resonance component 
$\omega^{(j)}=\omega^{(j)}_0\cdot n^{(j)}$
of the frequency system~$\omega$.

The classical version of (A.1) is the well-known classical
averaging method \cite{kar14,kar15}.
Namely, in our operator scheme:

--- the transformation $\wh{U}^{-1}_\ve\hat{f}\wh{U}_\ve$ 
is replaced by $\gamma^*_\ve f$,
where 
$$
\gamma^*_\ve=\exp(\ve\ad(f_0)+\ve^2\ad(f_1)+O(\ve^3))
$$ 
and $\ad(\cdot)$ denotes the Hamiltonian vector field;

--- the homological equations (A.4) and Eqs.~(A.2) are replaced
by 
\begin{equation} 
\begin{aligned}
\{H_0,f_0\}&=H_1-\underline{H_1},&\qquad
\{H_0,f_1\}&=H_2-\underline{H_2},
\\
\{H_0,\underline{H_1}\}&=0,&\qquad
\{H_0,\underline{H_2}\}&=0,
\end{aligned}
\tag{A.8}
\end{equation}
and the solution of (A.8) is given by (A.6), (A.7),
where one has to remove the sign $\,\widehat{}\,$ 
over all the functions;

--- the classical analog of (A.1) is 
$$
\gamma^*_{\ve} (\ad(H_0)+\ve\ad(H_1))
=(\ad(H_0)+\ve\ad(\underline{H_1})+\ve^2\ad(\underline{H_2})
+ O(\ve^3))\gamma^*_\ve,
$$
which means that by the symplectic transformation
$\gamma^{-1}_\ve$ one transforms the Hamiltonian field
corresponding to $H_0+\ve H_1$ into the Hamiltonian field
corresponding to $H_0+\ve\underline{H_1}+\ve^2\underline{H_2}$ 
with $\underline{H_1}$, $\underline{H_2}$ being in involution
with $H_0$. 

The flows of these Hamiltonians are related to each
other as follows: 
$$
\gamma^t_{H_0+\ve H_1}
=\gamma_\ve\circ\gamma^t_{H_0}\circ
\gamma^{\ve t}_{\underline{H_1}+\ve\underline{H_2}}\circ
\gamma^{-1}_\ve +O(\ve^3).
$$
This operator averaging scheme and its generalizations can be
found in \cite{n12}; see also more details and applications in
\cite{kar20}.

\end{document}